\documentclass[12pt]{article}
    \usepackage{srcltx}
   \usepackage{enumerate}
 \usepackage{latexsym}
   \usepackage{indentfirst}
\usepackage{amsmath,amsfonts,amssymb,amsthm}
   \usepackage[russian]{babel}

\newtheorem{theor}{Теорема}
\newtheorem{cor}{Следствие}
\newtheorem{ex}{ Пример}

\newtheorem{lem}{Лемма}

\newtheorem{predl}{П~р~е~д~л~о~ж~е~н~и~е}

\theoremstyle{remark}
\newtheorem{rem}
{З~а~м~е~ч~а~н~и~е~}

\title{
Asymptotic properties of one-step weighted $M$-estimators 
and applications to some regression problems
 }
\author{Yu. Yu. Linke
\thanks{Sobolev Institute of Mathematics,
4 Academician Koptyug ave., Novosibirsk, 630090  Russia;
Novosibirsk State University, 2 Pirogov str., 630090  Russia.
E-mail:
linke@math.nsc.ru}   }
\date{}

\begin{document}

 \maketitle

{\it Summary:}
 We study asymptotic behavior of one-step weighted $M$-estimators based on samples from arrays of not necessarily identically distributed random variables and representing
 explicit approximations to the corresponding consistent weighted $M$-estimators.  Sufficient conditions are presented for asymptotic normality of the one-step weighted $M$-estimators under consideration.
 As a consequence, we consider some well-known   nonlinear regression models where the procedure mentioned allow us to construct explicit asymptotically optimal estimators.

  %
%

%
%

\section{Введение}

{\bf 1.1.} Пусть  $X_1,\ldots, X_n$ --- независимые не обязательно одинаково распределенные наблюдения произвольной природы, распределения которых зависят от неизвестного параметра $\theta\in \Theta$.
Задача состоит в оценивании этого  параметра по  наблюдениям $X_1,\ldots, X_n$.
 Одним  из  основных общих методов  получения оценок 
  принято считать $M$-оценивание  (см., например, [\ref{2007-B}]). 
  Определим   $M$-оценки~$\widetilde\theta_n$  как
 статистики $\widetilde\theta_n$, которые определены с вероятностью, стремящейся к $1$, и являются решениями уравнений вида
\begin{equation}\label{p0-0+s}
\sum\limits_{i=1}^n  M_i(t,X_i)=0
\end{equation}
для некоторого набора функций   $  M_i(t,x)$, $i=1,\ldots,n,$ с условием  ${\bf E} M_i(\theta,X_i)=0$ при всех $i$.
Подчеркнем, что обычно функции $\{ M_i(t,x)\}$ в той или иной статистической постановке  выбираются таким образом, чтобы обеспечить желаемые свойства (в том или ином смысле оптимальность) соответствующей $M$-оценки 
   (см., например, [\ref{1997-H}], а также примеры в \S3).
\begin{rem}
Отметим, что далеко не все корни уравнения (\ref{p0-0+s}) можно определить на множестве выборок асимптотически полной меры. Так, например, возможна ситуация, когда количество корней при любом $n$ (и даже в пределе при $n\to\infty$, см. [\ref{1985-R}]) будет невырожденной случайной величиной, и лишь один корень можно определить с вероятностью, стремящейся к $1$. Более того,  возможны ситуации, когда несколько корней $M$-уравнения могут быть определены с вероятностью, стремящейся к $1$, или таковых корней нет, но с ненулевой вероятностью другие корни могут  существовать (см. подробности в [\ref{2015-2?}]).
\end{rem}

 Хорошо известно, что
     вычисление $M$-оценок  
     нередко связано с серьезными техническими трудностями (например, вычисление  оценок максимального правдоподобия или  оценок  метода наименьших квадратов для задач нелинейной регрессии), особенно в случае существования большого числа корней уравнения (\ref{p0-0+s}). 
  Ситуация существенно упрощается, если известна некоторая  предварительная состоятельная  оценка~$\theta_n^*$ параметра $\theta$, приближающая параметр с нужной нам скоростью сходимости.
    Оказывается, что в этом случае достаточно лишь одного шага  метода Ньютона, 
    чтобы в явном виде построить оценку, обладающую  теми  же асимптотическими  свойствами, как и  $M$-оценка. Так называемую {\it одношаговую $M$-оценку}   определим соотношением
\begin{equation}\label{p0-0s}
\theta_{n}^{**}
=\theta_n^*-{\sum\limits_{i=1}^n  M_i(\theta_n^*,X_i)}\big/{\sum\limits_{i=1}^n  M_i'(\theta_n^*,X_i)},
\end{equation}
где $\theta_n^*$  (здесь и далее) --- некоторая предварительная состоятельная оценка параметра $\theta$, а через $ M_i'(t,X_i)$ обозначены производные этих функций по первому аргументу.

Оценка  $\theta_{n}^{**}$ представляет собой  одну итерацию метода Ньютона   с начальной точкой $t=\theta_n^*$ для приближенного вычисления корня
уравнения  (\ref{p0-0+s}). В силу состоятельности $\theta_n^*$ этот корень в широких условиях  будет ближайшим корнем к $\theta$ (и этот корень, тем не менее, может не быть $M$-оценкой; см. [\ref{2015-2?}], а также замечание~\ref{p1-usl2}).
Отметим, что в литературе на протяжении нескольких последних  десятилетий, начиная с работ  П. Хьюбера,
П. Бикела, Р. Серфлинга, Н. Веравербеке,  П.-К. Сена, К. Мюллер,   Х. Юречковой,  Д. Рупперта, Р. Кэрролла и др. (см. библиографию здесь и 
 в [\ref{2015-2?}])
для обозначения таких по  сути  {\it двухшаговых} оценок
  используется термин {\it одношаговые} оценки в том смысле, что
при наличии некоторой {\it предварительной} оценки (существование которой зачастую лишь постулируется) именно  {\it за один шаг}
 метода Ньютона или некоторых его  модификаций
происходит качественное  улучшение этой начальной статистики.
Следуя устоявшейся терминологии оценки такого типа будем называть {\it одношаговыми}. 

При некоторых условиях регулярности (см., например, [\ref{2007-B}] и [\ref{2015-2?}]) 
 для $M$-оценки $\widetilde\theta_n$ из (\ref{p0-0+s}) и одношаговой $M$-оценки $\theta_{n}^{**},$ определенной в
 (\ref{p0-0s}),  имеют место сходимости
\begin{equation}\label{p0-0-s}
\begin{split}
\frac{J_n}{\sqrt{I_n}}\big(\widetilde\theta_{n}-\theta\big)\Longrightarrow {\cal N}(0,1),\qquad
\frac{J_n}{\sqrt{I_n}}\big(\theta_{n}^{**}-\theta\big)\Longrightarrow {\cal N}(0,1),\quad \mbox{где}\\
I_n:=\sum\limits_{i=1}^n{\bf E} M_i^2(\theta,X_i)>0, \qquad J_n:=\sum\limits_{i=1}^n{\bf E} M_i'(\theta,X_i)\neq 0,
\end{split}
\end{equation}
 а запись вида
$\eta_n\Rightarrow {\cal N}(0,1)$ означает слабую
сходимость  распределений $\eta_n$  к стандартному  нормальному закону (всюду пределы, если не оговорено противное, берутся при $n\to \infty $).
Таким образом, одношаговая $M$-оценка $\theta_n^{**}$ оказывается  асимптотически нормальной с той же асимптотической дисперсией, что и $M$-оценка $\widetilde\theta_n$. При этом  в (\ref{p0-0-s}) соотношение для $\theta_n^{**}$ может иметь место, в то время как $M$-оценка $\widetilde\theta_n$ не существует или не состоятельна  (впервые подобный эффект был отмечен  Л.Ле Камом [\ref{1956-L}] в случае однородной выборки,  подробности и комментарии в случае не обязательно одинаково распределенной выборки приведены в [\ref{2015-2?}]).

{\bf 1.2. } Предположим теперь, что функции $M_i(t,X_i)$ из (\ref{p0-0+s}) имеют вид  
%
%
\begin{equation}  \label{p1-170+ss}
M_i(t,X_i)=h_i(t) \widetilde  M_i(t,X_i).
 \end{equation}
 Подобная факторизация для функций $ M_i(t,X_i)$ естественным образом возникает в некоторых статистических задачах. Например, оценка метода наименьших квадратов в задаче нелинейной регрессии есть $M$-оценка $\widetilde\theta_n$, 
 являющаяся решением уравнения
 $\sum\limits_{i=1}^nf_i'(t)\big(X_i-f_i(t)\big)=0, $ где $f_i(\cdot)$ ---  известная регрессионная функция (см. подробности и др. примеры в \S3). Тем самым, можно положить $h_i(t)=f_i'(t)$ и $\widetilde M_i(t,X_i)=X_i-f_i(t)$.

Если функции $h_i(t)$ достаточно регулярные (скажем, удовлетворяют условию Гельдера), то в окрестности истинного значения параметра $\theta$ функции
$ M_i(t,X_i)$ вида (\ref{p1-170+ss}) и  $h_i(\theta_n^*) \widetilde M_i(t,X_i)$ будут в известном смысле близкими. Это позволяет надеяться, что в окрестности истинного значения параметра   будут близки  корни $\widetilde \theta_n$
уравнения (\ref{p0-0+s}) в случае (\ref{p1-170+ss}) и корни $\widehat\theta_n$ уравнения
 \begin{equation} \label{p1-10+}
\sum\limits_{i=1}^n h_i(\theta_n^*) \widetilde M_i(t,X_i)=0.
\end{equation}

Используя     один шаг   метода Ньютона   с начальной точкой $t=\theta_n^*$ для приближенного поиска состоятельного решения $\widehat\theta_n$ уравнения (\ref{p1-10+}),
получаем следующую одношаговую оценку для параметра $\theta$:
 \begin{equation}  \label{p1-172}
 \widetilde\theta_{n}^{**}=\theta_n^*- \sum\limits_{i=1}^n
h_i(\theta_n^*) \widetilde M_{i}(\theta_n^*,X_i)\Big/\sum\limits_{i=1}^nh_i(\theta_n^*)\widetilde  M_{i}'(\theta_n^*,X_i).
 \end{equation}
  Корни $\widehat\theta_n $ уравнения вида (\ref{p1-10+}), определенные с вероятностью, стремящейся к~$1$, будем называть    {\it   взвешенными  $M$-оценками}, а оценки, определяемые соотношениями вида   (\ref{p1-172}) ---  {\it  одношаговыми взвешенными  $M$-оценками}.

%
%
%
Отметим, что в случае (\ref{p1-170+ss}) одношаговая $M$-оценка $\theta_n^{**}$ из (\ref{p0-0s}) имеет  вид
 \begin{equation}  \label{p1-171}
 \theta_{n}^{**}=\theta_n^*- \frac{\sum\limits_{i=1}^n
h_i(\theta_n^*)\widetilde M_{i}(\theta_n^*,X_i)}{\sum\limits_{i=1}^nh_i(\theta_n^*) \widetilde M_{i}'(\theta_n^*,X_i)+
\sum\limits_{i=1}^nh_i'(\theta_n^*) \widetilde M_{i}(\theta_n^*,X_i)}.
 \end{equation}
При выполнении некоторых  условий регулярности  эта одношаговая $M$-оценка
 является асимптотически нормальной с асимптотической дисперсией, определяемой соотношением (см. (\ref{p0-0-s}) и (\ref{p1-170+ss}))
 \begin{equation}  \label{p1-2323}
 {I_n}/{J_n^2}\equiv
{\sum\limits_{i=1}^nh_i^2(\theta){\bf E} \widetilde M_{i}^2(\theta,X_i)}\Big/{\Big(\sum\limits_{i=1}^nh_i(\theta){\bf E} \widetilde M_{i}'(\theta,X_i)\Big)^2},
 \end{equation}
где мы учли, что  ${\bf E} \widetilde M_{i}(\theta,X_i)=0$.

 В данной работе изучаются одношаговые взвешенные $M$-оценки. Будет установлено, что при некоторых условиях регулярности одношаговая взвешенная $M$-оценка 
 из
(\ref{p1-172}) является асимптотически нормальной с той же асимптотической дисперсией (\ref{p1-2323}), что и оценка $\theta_n^{**}$ из (\ref{p1-171}).
 Отметим, что одношаговая взвешенная $M$-оценка 
 из (\ref{p1-172}) проще, и ее конструкция  не предполагает, в отличие от оценки  $\theta_{n}^{**}$ из (\ref{p1-171}), дифференцируемость функций $h_i(t)$.  
В частности,  для отмеченной выше задачи  одношагового приближения  оценки метода наименьших квадратов последнее  означает меньшие условия на гладкость регрессионной функции $f_i(\theta)$.
  Другие  аргументы в пользу рассмотрения одношаговых взвешенных $M$-оценок мы приведем в \S3.

{\bf 1.3.} Идея  одношагового оценивания впервые  была предложена  Р. Фишером [\ref{1925-F}] в задаче приближенного поиска состоятельных оценок максимального правдоподобия в случае однородной выборки. Исследование одношаговых оценок Фишера, а также их уточнений и   обобщений в случае одинаково распределенной выборки содержится в ряде монографий и статей (см., например, [\ref{2007-B}],  [\ref{1975-Z}], [\ref{1985-J}], [\ref{1990-J}], [\ref{1991-L}], [\ref{1984-H}], [\ref{1956-L}], [\ref{1980-S}], [\ref{2007-V}], [\ref{2014-1}], [\ref{2015-1?}]  и ссылки там же).

Возможность использования методологии одношагового оценивания в статистических задачах с разнораспределенными выборочными наблюдениями хорошо  известна (см., например, монографии [\ref{2007-B}], [\ref{1984-H}], [\ref{2006-J}], [\ref{2012-J1}], [\ref{1997-H}]) 
и успешно  используется  в различных специальных статистических задачах
  (см., например,  [\ref{1975-B}], [\ref{1984-H}],
[\ref{1988-G}],
 [\ref{1987-J}] 
 [\ref{1999-F1}], [\ref{2011-B}],
  [\ref{1999-F}], [\ref{2006-J}],
  [\ref{2012-J1}], [\ref{1994-M1}] 
  и ссылки там же).
  По-видимому, впервые одношаговые оценки в случае разнораспределенных наблюдений были предложены в
[\ref{1975-B}] для задачи робастного оценивания  в линейных моделях.
В работе автора   [\ref{2015-2?}] 
изучается асимптотическое поведение одношаговых $M$-оценок общего вида (функции, определяющие эти оценки, не связаны со специальной статистической задачей).  
Терминология одношаговых взвешенных $M$-оценок и взвешенных $M$-оценок  использовалась, например, в [\ref{2011-B}] для схожих  статистик в  специальной обобщенной модели линейной регрессии, а также частично  в [\ref{1988-G}] для одной частной регрессионной задачи.

  Как правило, в литературе исследуются  два аспекта, связанные с одношаговыми оценками:  сближение  одношаговых  оценок и параметра (в частности, асимптотическая нормальность этих оценок), а также  сближение одношаговой $M$-оценки и состоятельной $M$-оценки.
Стоит отметить,  что обычно при 
изучении одношаговых оценок   в качестве предварительных оценок параметра рассматриваются $\sqrt{n}$-ограниченные оценки (т.е.   оценки, удовлетворяющих соотношению
$\sqrt{n}(\theta_n^*-\theta)=O_p(1)$).
Исключение составляют лишь монография  [\ref{1975-Z}] и работы [\ref{2014-1}], [\ref{1987-J}], [\ref{2015-2?}], [\ref{2015-1?}].
 Тем не менее, задачи приближенного поиска  $M$-оценок возникают и в случае, когда скорость сближения предварительной оценки и  параметра оказывается  медленнее, чем $O(n^{-1/2})$ (см., например, [\ref{2015-1?}], [\ref{2014-1}]). В работе автора  [\ref{2015-2?}]  найдены условия асимптотической нормальности одношаговых $M$-оценок общего вида  в широком спектре ограничений на точность предварительной оценки.
%
%
Более подробный  библиографический обзор об  одношаговом оценивании   содержится в [\ref{2015-2?}]. 

{\bf 1.4.} О структуре работы.  В \S 2  при достаточно общих условиях (в том числе и  на точность предварительной оценки $\theta_n^*$)
 доказана асимптотическая нормальность одношаговых взвешенных $M$-оценок.
Построены оценки для асимптотической дисперсии одношаговых взвешенных  $M$-оценок    и доказана сходимость к стандартному нормальному закону погрешностей оценивания при замене асимптотических дисперсий их оценками, 
 что может быть полезным при построении доверительных интервалов и для проверки статистических  гипотез.
Проведена процедура оптимизации функций $\{h_i(\cdot)\}$, участвующих в построении одношаговых взвешенных $M$-оценок (а также взвешенных $M$-оценок).

  В дальнейшем символ <<тильда>> у функций $\widetilde M_i(t,X_i)$ в определении   (\ref{p1-172}) одношаговых взвешенных $M$-оценок   
  мы будем опускать, так что формально одношаговые взвешенные $M$-оценки 
  можно рассматривать как  некоторые обобщения классических одношаговых $M$-оценок (\ref{p0-0s}).
 В разделе 2.3 в качестве следствия 
  приведены условия асимптотической нормальности одношаговых $M$-оценок. 
Приведенные результаты в  случае однородной выборки 
обобщают и    утверждения из [\ref{1980-S}], [\ref{1975-Z}], [\ref{1991-L}] и [\ref{2014-1}]  об асимптотической нормальности оценок Фишера и одношаговых $M$-оценок, построенных по выборке  одинаково распределенных наблюдений (см. подробности в [\ref{2014-1}] и [\ref{2015-2?}]).

%

 \S3 целиком посвящен регрессионным приложениям.   Наш интерес к одношаговому оцениванию в случае разнораспределенных наблюдений в первую очередь связан с возможностью эффективного применения этих процедур в различных задачах оптимального оценивания параметров  в моделях нелинейной регрессии. В частности, при  наличии  некоторой явной предварительной оценки параметра одношаговые   оценки позволяют находить явные  асимптотически оптимальные оценки неизвестного параметра в  гетероскедастических моделях нелинейной регрессии. 
Предлагаемая процедура оценивания иллюстрируется  несколькими  примерами   известных  задач нелинейной регрессии. 

 Доказательства всех утверждений отнесены в~\S4.

\section{Основные результаты.}

{\bf 2.1.} В этом разделе будут получены условия асимптотической нормальности   одношаговых взвешенных  $M$-оценок.

Нам потребуются следующие предположения. 

$ ({\bf R}_1).$
Наблюдается выборка объема $n$,  состоящая из независимых элементов   $X_{1},$ $X_{2},\ldots, X_{n}$  со значениями в
 произвольном измеримом пространстве
${\cal X}$ и  распределениями~${\cal L}_1,$ ${\cal L}_2,\ldots,{\cal L}_n$, зависящими  от    интересующего нас  основного
неизвестного параметра $\theta\in\Theta\subset
\mathbb{R}$, где $\Theta$ --- открытое множество.  Кроме того,  эти распределения зависят, вообще говоря, от $n$, и, возможно, зависят  также  от некоторых мешающих  параметров
 произвольной природы.
%
%

$ ({\bf R}_2).$  На множестве  $\Theta$ при любом  $i$ заданы зависящие, вообще говоря, от $n$  функции $h_i(t)$, $M_i(t,X_i)$ и $M'_i(t,X_i)$, при этом для каждого интервала
$(t_1,t_2)$, целиком лежащего в $\Theta$, с вероятностью $1$ имеет место равенство
\begin{equation}\label{p1-s1}
M_i(t_2,X_i)-
M_i(t_1,X_i)=\int\limits_{t_1}^{t_2}
 M'_i(t,X_i)dt \qquad \forall (t_1,t_2)\subset  \Theta,
 \end{equation}
  справедливы  соотношения
\begin{equation}\label{p1-25}
{\bf E}M_i(\theta,X_i)=0, \qquad {\bf E}M^2_i(\theta,X_i)<\infty, \qquad
{\bf E}|M'_i(\theta,X_i)|<\infty,
\end{equation}
а  функции  $h_i(\cdot)$ 
 удовлетворяют  условию Гельдера с показателем $p\in(0,1]$: 
\begin{equation}\label{p1-j41w}
|h_i(t_1)-h_i(t_2)|\leq \overline h_i|t_1-t_2|^p \qquad \forall t_1,t_2\in\Theta.
\end{equation}

%

%
%

%
%
%

$ ({\bf R}_3).$
   Начиная с некоторого $n$
\begin{equation}\label{p1-u25}
I_{n,h}:=\sum\limits_{i=1}^n
h_i^2(\theta){\bf E}M^2_i(\theta,X_i)>0, \qquad J_{n,h}:=
\sum\limits_{i=1}^n h_i(\theta){\bf E}M_i'(\theta,X_i)\neq 0
\end{equation}
и имеют место сходимости: $\quad |J_{n,h}|/\sqrt{I_{n,h}}\to \infty, $
\begin{equation}\label{p1-5-}
\sum\limits_{i=1}^nh_i(\theta) M'_i(\theta,X_i)/J_{n,h}\stackrel{p}{\to}1, \qquad
\sum\limits_{i=1}^n h_i(\theta)M_i(\theta,X_i)/\sqrt{I_{n,h}}\Longrightarrow {\cal N}(0,1).
\end{equation}

%
%

Положим $
\overline\tau_i(\delta):={\bf E}\tau_i(\delta,X_i),
$
  где
\begin{equation}   \label{p1-tau}
\tau_{i}(\delta,X_i):=
\begin{cases}
\sup\limits_{t
:\; |t-\theta|\leq\delta}\big|
M'_i(t,X_i)-M'_i(\theta,X_i) \big|,\quad \mbox{если }\; [\theta-\delta,\theta+\delta]\subset\Theta,
\\
\; \infty,\quad \mbox{иначе}.
\end{cases}
\end{equation}

$ ({\bf R}_4).$   Имеет место сходимость
$$\lim\limits_{\delta\to 0}\limsup\limits_{n\to \infty}\sum\limits_{i=1}^n\big(|h_i(\theta)|+\delta^p\overline h_i\big) \overline\tau_{i}(\delta)/|J_{n,h}|= 0. 
$$

$ ({\bf R}_5).$ Оценка  $\theta_n^*$ такова, что 
 \begin{gather}                                                            
\label{p1-27!!}
|\theta_n^*-\theta|\sum\limits_{i=1}^n\big(\overline h_i|\theta_n^*-\theta|^p+|h_i(\theta)|\big)\overline\tau_i(|\theta_n^*-\theta|)\big/\sqrt{I_{n,h}}\stackrel{p}{\to}  0,\\
                                                                                     \label{p1-j27}
|\theta_n^*-\theta|^{p}\sum\limits_{i=1}^n\overline h_i{\bf E}|M_i'(\theta,X_i)|\big/|J_{n,h}|\stackrel{p}{\to}  0.
  \end{gather}

Обозначим через $(R)$ группу условий, которая  включает в себя   предположения $(R_1),\ldots, (R_5)$.
Одношаговые взвешенные $M$-оценки определим соотношением
\begin{gather} \label{p1-7+}
 \theta_{n,M,h}^{**}=\theta_n^*- {\sum\limits_{i=1}^nh_i(\theta_n^*) M_i(\theta_n^*,X_i)}\big/{
 \sum\limits_{i=1}^nh_i(\theta_n^*)M'_i(\theta_n^*,X_i)}.
 \end{gather}

Имеет место
\begin{theor}                                     \label{p1-t1-2+}
Пусть  выполнено  предположение $(R)$ и условие
\begin{equation}  \label{p1-j41}
\delta_{n}^h:=
\sum\limits_{i=1}^n\big(h_i(\theta_n^*)-h_i(\theta)\big)M_i(\theta,X_i)/\sqrt{I_{n,h}}\stackrel{p}{\to}0.
\end{equation}
 Тогда 
  одношаговая взвешенная  $M$-оценка $\theta_{n,M,h}^{**}$  определена с вероятностью, стремящейся к $1$, и
\begin{gather}                                                                                           \label{p1-24-}
\frac{J_{n,h}}{\sqrt{I_{n,h}}}\big(\theta_{n,M,h}^{**}-\theta\big)\Longrightarrow {\cal N}(0,1).
 \end{gather}
%
\end{theor}

\begin{rem}\label{p1-usl1}
Достаточные условия для  выполнения соотношений  (\ref{p1-5-}), являющихся вариантами закона больших чисел и центральной предельной теоремы в схеме серий, хорошо известны  (см., например, [\ref{2007-B}]).
Для справедливости  $(R_4)$ и $(R_5)$ достаточно, чтобы имели место следующие условия:
   $$
|M_i(t,X_i)-M_i(\theta,X_i)|\leq \overline{M_i'}|t-\theta|^q \quad \mbox{при}\; q\in(0,1],\qquad
n^{\frac{1}{ 2(1+q)}}(\theta_n^*-\theta)\stackrel{p}{\to}0,
$$
\begin{equation}\label{p1-50w}
1/I_{n,h}+1/|J_{n,h}|=O\big(1/n\big),\quad \sup\limits_{n,i}\{|h_i(\theta)|,\overline h_i, \;{\bf E}\overline{M_i'},\; {\bf E}|M_i'(\theta,X_i)| \}<\infty, \; p>0.
\end{equation}
\end{rem}

\begin{rem}
Обсудим подробно  условие (\ref{p1-j41}), используемое в теореме \ref{p1-t1-2+}.
 Если  функции $h_i(t)$ достаточно гладкие, чтобы при выводе
 (\ref{p1-j41})  воспользоваться формулой Тейлора, то не трудно получить простые  достаточные условия для~(\ref{p1-j41}).    В случае же, когда функции $h_i(t)$ удовлетворяют лишь условию Гельдера,  доказательство  (\ref{p1-j41}) представляет некоторую техническую трудность. Приведем два утверждения, содержащие  достаточные условия для выполнения  (\ref{p1-j41}).
 В предложениях \ref{p1-pre1} и  \ref{p1-pre2}  считаем, что  выполнено условие $(R_1)$,  первые два  соотношения в  {\rm(\ref{p1-25})} и первое условие в (\ref{p1-u25}). Без ограничения общности будем также  считать, что функции $h_i(t)$
продолжены на $\mathbb{R}$ с выполнением  условия Гельдера~(\ref{p1-j41w}) на всей числовой прямой.

 Из теоремы 3 в [\ref{2011-1}] нетрудно извлечь
 \begin{predl}\label{p1-pre1}
 Пусть  оценка $\theta_n^*$
имеет следующую структуру
\begin{equation}\label{p1-28w}
\theta_n^*-\theta={\sum\limits_{i=1}^n u_{ni}(\theta,X_i)}\Big/{\Big(1+\sum\limits_{i=1}^n v_{ni}(\theta,X_i)\Big)},
\end{equation}
при этом ${\bf E}u_{ni}(\theta,X_i)= {\bf E}v_{ni}(\theta,X_i)=0$, ${\bf E}u_{ni}^2(\theta,X_i)<\infty$, $ {\bf E}v_{ni}^2(\theta,X_i)<\infty$ и   
\begin{eqnarray}\label{p1-29w}
 d_{n} +\sum\limits_{i=1}^n {\bf E}v_{ni}^2(\theta,X_i)\to 0\quad \mbox{при}\quad d_n^2:=\sum\limits_{i=1}^n {\bf E}u_{ni}^2(\theta,X_i),\nonumber \\
d_n^{2p}\left(\sum\limits_{i=1}^n \Big(\overline h_i^2{\bf E} M_{i}^2(\theta,X_i)\Big)^{1/(2-p)}\right)^{2-p}\Big/I_{n,h}\to 0.
 \end{eqnarray}
Тогда имеет место {\rm(\ref{p1-j41})}.
 \end{predl}
 Подчеркнем, что во всех примерах  \S3 оценки
 $\theta_n^*$ имеют структуру (\ref{p1-28w}).
 Отметим еще, что для справедливости центрального условия (\ref{p1-29w}) достаточно, чтобы выполнялись соотношения
\begin{gather}\label{p1-31w}
 n^{(1/p-1)/2}d_n\to 0,\quad {1}/I_{n,h}=O\big(1/n\big),\quad \sup\limits_{n,i}\{\overline h_i,\; {\bf E} M_i^2(\theta,X_i)\}<\infty.
 \end{gather}
В частности, первое условие в (\ref{p1-31w}) выполнено в следующих трех случаях:
$$
\sqrt{n}d_n=O(1)\; \mbox{и}\; p>1/2, \qquad n^{1/4}d_n\to 0\; \mbox{и}\; p\geq 2/3,\qquad d_n\to 0\; \mbox{и}\; p=1.
$$
В случае, когда $\theta_n^*$ имеет  дробно-линейную структуру,  достаточные условия для (\ref{p1-j41}) можно получать и в терминах
сближения $\theta_n^*$ и $\theta$ по вероятности (т.е. в терминах условия  $n^{(1/p-1)/2}(\theta_n^*-\theta)\stackrel{p}{\to} 0$) (см. подробности в [\ref{2013-2}]).

 Обозначим через~ ${\bf E}_{\neq i}$    условное
математическое ожидание,
  взятое при условии, что 
  при всех $j\neq i$ фиксированы значения
    $X_{j},$ $j=1,\ldots,n$.
В следующем предложении  достаточные условия для   (\ref{p1-j41}) не содержат  предположения о дробно-линейной структуре $\theta_n^*.$
 \begin{predl}\label{p1-pre2}
 Если  ${\bf E}|\theta_n^*|^2<\infty $ и
\begin{gather}\label{p1-j28}
\big({\bf E}|\theta_n^*-\theta|^2\big)^p\sum\limits_{i=1}^n\overline h_i^2{\bf E}M_i^2(\theta,X_i)\big/I_{n,h}\to 0,
\\ \label{p1-j29}
\sum\limits_{i=1}^n\overline h_i\big({\bf E}|\theta_n^*-{\bf E}_{\neq i}\theta_{n}^*|^2\big)^{p/2}\big({\bf E}M_i^2(\theta,X_i)\big)^{1/2}\big/\sqrt{I_{n,h}}\to 0,
  \end{gather}
  то имеет место {\rm(\ref{p1-j41})}.~\hfill$\square$
 \end{predl}%

\end{rem}

Положим
\begin{equation} \label{p1-9-s}
 d_{n,M,h}^*={\sum\limits_{i=1}^n h_i(\theta_n^*)M_i'(\theta_n^*,X_i)}\Big/{
\Big(\sum\limits_{i=1}^n h_i^2(\theta_{n,M,h}^{**})M_i^2(\theta_{n,M,h}^{**},X_i)\Big)^{1/2}}
\end{equation}
 Следующее утверждение  может быть полезным при построении доверительных интервалов и проверке гипотез, поскольку  множитель  $d_{n,M,h}^*$, нормирующий  разность   $\theta^{**}_{n,M,h}-\theta$, является статистикой.

\begin{theor}                                     \label{p1-t21}
Если   выполнены  условия теоремы~{\rm\ref{p1-t1-2+}} и
\begin{equation} \label{p1-8-}
\frac{\sum\limits_{i=1}^n h_i^2(\theta)M_i^2(\theta,X_i)}{I_{n,h}}\stackrel{p}{\to}1,\quad
\frac{\sum\limits_{i=1}^n \big(\overline h_i^2+h_i^2(\theta)\big)[M_i'(\theta,X_i)]^2}{J_{n,h}^2}
\stackrel{p}{\to}0, \quad \frac{\sum\limits_{i=1}^n \overline h_i^2M_i^2(\theta,X_i)}{J_{n,h}^{2p}I_{n,h}^{1-p}}
\stackrel{p}{\to}0,
\end{equation}
то статистика $d_{n,M,h}^*$
 определена с вероятностью, стремящейся к $1$ и
\begin{equation} \label{p1-9-ss}
d_{n,M,h}^*\big(\theta_{n,M,h}^{**}-\theta\big)\Longrightarrow {\cal N}(0,1).
\end{equation}
\end{theor}
 Дополнительные условия в теореме \ref{p1-t21} не являются ограничительными, поскольку первая  сходимость в (\ref{p1-8-}) --- это вариант закона больших чисел, а для справедливости двух других сходимостей из  (\ref{p1-8-}) достаточно, чтобы вместе с первым соотношением в   (\ref{p1-50w}) выполнялось условие  $
\sup\limits_{n,i}\{|h_i(\theta)|, \overline h_i,\; {\bf E}|M_i'(\theta,X_i)|,\;{\bf E}M_i^2(\theta,X_i) \}<\infty.
$

{\bf 2.2.} Рассмотрим теперь вопрос оптимизации функций $h_i(\cdot),$  напрямую связанный с
минимизацией асимптотической дисперсии $I_{n,h}/J_{n,h}^2.$

\begin{predl}\label{p1-t4-5}
Пусть $J_{n,h}\equiv\sum\limits_{i=1}^n h_i(\theta){\bf E}M_i'(\theta,X_i)\neq 0$ начиная с
некоторого $n$. В этом
случае определена величина ${I_{n,h}}/{J_{n,h}^2}$ и 
\begin{equation}\label{p1-391}
\frac{I_{n,h}}{J_{n,h}^2}\geq \frac{I_{n,h^o}}{J_{n,h^o}^2}
\equiv \left(\sum\limits_{i=1}^n
\frac{\big({\bf E}M_i'(\theta,X_i)\big)^2}{{\bf E}M_i^2(\theta,X_i)}\right)^{-1}\quad\mbox{при}\quad
h_i^{o}(\theta)= \frac{{\bf E}M_i'(\theta,X_i)}{{\bf E}M_i^2(\theta,X_i)}.
\end{equation}%
Если дополнительно
${\rm
sign}h_i(\theta)={\rm sign}h_i^{o}(\theta)$, то
справедливо неравенство:
\begin{equation}                         \label{p1-393}
1\leq\frac{I_{n,h}/J_{n,h}^2}
{I_{n,h^o}/J_{n,h^o}^2} \leq
1+\frac{(\sqrt{H/h}-1)^2}{2\sqrt{H/h}} \le\sqrt{\frac{H}{h}},\quad\mbox{где}
\end{equation}
 \begin{gather}\label{p1-394}
 h:=\min\limits_{\{i:\; {\bf E}M_{i}'(\theta,X_i)\neq0\}}
h_{i}(\theta)\big/h_{i}^{o}(\theta)>0,
\quad H:=\max\limits_{\{i:\; {\bf E}M_{i}'(\theta,X_i)\neq0\}}
h_{i}(\theta)\big/ h_{i}^{o}(\theta)<\infty.
\end{gather}
\end{predl}
Неравенство (\ref{p1-393}) можно интерпретировать как некоторое
свойство устойчивости оценок с асимптотической дисперсией $I_{n,h}/J_{n,h}^2$
как функционалов,
зависящих от набора  функций $\{h_{i}(\cdot)\}$: чем лучше
$h_{i}(\theta)$ приближают оптимальные
$h_{i}^{o}(\theta),$ тем меньше асимптотическая дисперсия
$I_{n,h}/J_{n,h}^2$  отличается от
оптимальной $I_{n,h^o}/J_{n,h^o}^2$.
Ранее аналог  неравенства (\ref{p1-393})
для асимптотических дисперсий специального вида был получен в [\ref{2008}] и [\ref{2009}].

Неравенство  (\ref{p1-393})  позволяет сделать
 некоторые  рекомендации по выбору функций $h_i(\cdot),$
участвующих в определении  оценки $\theta_{n,M,h}^{**}.$  Если оптимальные  функции  $ h_{i}^o(\cdot) $  нам известны (см., например, \S3), то в качестве  $ h_{i}(\cdot)  $ разумно взять
оптимальные функции, т.е. положить  $h_i(\cdot)= h_{i}^{o}(\cdot)$.
Но  если, скажем, оптимальные функции $h_{i}^{o}(\cdot)  $
 не удовлетворяют  условиям теоремы    \ref{p1-t1-2+}, 
 то можно использовать функции $
\widetilde h_{i}^{o}(\cdot)$,  полученные в результате
соответствующего сглаживания функций $ h_{i}^{o}(\cdot)$. Кроме того, поскольку распределения $\{{\cal L}_i\}$ могут зависеть от некоторых мешающих параметров, то  оптимальные  функции $  h_{i}^o(\cdot)  $  могут быть не известны (см., например, замечание \ref{p1-zaz}). В этом случае   можно рекомендовать
использовать некоторые функции $\widetilde{h}_{i}^o(\cdot)  $, которые  должны, по-возможности, ``не очень
сильно отличаться'' от неизвестных функций $  h_{i}^o(\cdot)$.
При этом соотношение (\ref{p1-393}) гарантирует, что асимптотическая дисперсия $I_{n,\widetilde h^o}/J_{n,\widetilde h^o}^2$ одношаговой взвешенной $M$-оценки будет
``не очень сильно отличаться'' от оптимальной асимптотической  дисперсии $I_{n, h^o}/J_{n, h^o}^2$.

Разумеется, приведенные здесь соображения по выбору  функций $h_i(\cdot)$ в полной мере относятся и к выбору этих функций при построении
$M$-оценок в случае (\ref{p1-170+ss}) и  взвешенных $M$-оценок.

{\bf 2.3.} В качестве следствия из теоремы   \ref{p1-t1-2+} можно  получить  условия асимптотической нормальности
для одношаговых $M$-оценок (\ref{p0-0s}). 
Чтобы привести соответствующее утверждение, нам потребуется аналог условия  $(R)$ при   $h_i(t)\equiv 1$ и $\overline h_i\equiv 0$.

Обозначим через $(R^-)$ группу условий, включающую  в себя
 условие $ ( R_1)$, условия (\ref{p1-s1}) и (\ref{p1-25}) из  $ ( R_2)$ и приводимые далее условия $ ( R_3^-)$-$ ( R_5^-)$.

$ ({\bf R}_3^-).$
   Начиная с некоторого $n$
$I_{n}:=\sum\limits_{i=1}^n
{\bf E}M^2_i(\theta,X_i)>0$, $J_{n}:=
\sum\limits_{i=1}^n {\bf E}M_i'(\theta,X_i)\neq 0$
и имеют место сходимости: $$
\quad |J_n|/\sqrt{I_n}\to \infty, \qquad
\sum\limits_{i=1}^n M'_i(\theta,X_i)/J_n\stackrel{p}{\to}1, \qquad
\sum\limits_{i=1}^n M_i(\theta,X_i)/\sqrt{I_n}\Longrightarrow {\cal N}(0,1).
$$


$ ({\bf R}_4^-).$   Имеет место соотношение
$\lim\limits_{\delta\to 0}\limsup\limits_{n\to \infty}\displaystyle\sum\limits_{i=1}^n \overline\tau_{i}(\delta)/|J_n|= 0. 
$

$ ({\bf R}_5^-).$ Оценка  $\theta_n^*$ такова, что
$
|\theta_n^*-\theta|\sum\limits_{i=1}^n\overline\tau_i(|\theta_n^*-\theta|)\big/\sqrt{I_n}\stackrel{p}{\to}  0.
 $

Из теоремы  \ref{p1-t1-2+} нетрудно извлечь

\begin{cor} \label{p1-cort1}
Пусть  выполнено  предположение $(R^-)$. Тогда
 одношаговая  $M$-оценка {\rm(\ref{p0-0s})}
%
определена с вероятностью, стремящейся к $1$, и $$
\frac{J_n}{\sqrt{I_n}}\big(\theta_{n}^{**}-\theta\big)\Longrightarrow {\cal N}(0,1).
 $$
  Если дополнительно
 $
\sum\limits_{i=1}^n M_i^2(\theta,X_i)\big/I_n\stackrel{p}{\to}1$ и $
\sum\limits_{i=1}^n [M_i'(\theta,X_i)]^2\big/J_n^2\stackrel{p}{\to}0,
$
то 
$$
d_{n}^*\big(\theta_{n,M}^{**}-\theta\big)\Longrightarrow {\cal N}(0,1)\quad\mbox{при}\quad d_{n}^*={\sum\limits_{i=1}^n M_i'(\theta_n^*,X_i)}\Big/{
\Big(\sum\limits_{i=1}^nM_i^2(\theta_{n,M}^{**},X_i)\Big)^{1/2}}.
 $$
 \end{cor}

\begin{rem}
  Условие  $(R_5)$ и его частный случай  $(R_5^-)$ являются  некоторыми универсальными ограничениями, связывающими гладкость функций,  
   определяющих одношаговую оценку, 
  со скоростью сближения предварительной оценки $\theta_n^*$  и параметра $\theta$,
    которые
нужны для  асимптотической нормальности  
одношаговых  оценок.  
При этом точность предварительной оценки $\theta_n^*$ и гладкость функций
$\{M_i(\cdot,X_i)\}$   обратно пропорциональны друг другу:  чем меньше точность, тем больше должна быть гладкость
 (см. также замечание \ref{p1-usl1}).

Например,  для выполнения $(R_5^-)$ достаточно, чтобы
 функции $\{M_i(t,X_i)\}$  и оценка $\theta_n^*$
одновременно удовлетворяли одному из следующих двух условий:
\begin{eqnarray*}
\limsup\limits_{n\to\infty}  \frac{\sum\limits_{i=1}^n\overline\tau_{i}(\delta)}{|J_n|}=o(\delta^\alpha)\;\mbox{при}\;\delta\to 0, \quad
\Big(\frac{|J_n|}{\sqrt{I_n}}\Big)^{1/(1+\alpha)}|\theta_n^*-\theta|=O_p(1),\quad 0\leq \alpha<1;
\\
\limsup\limits_{n\to\infty}  \frac{\sum\limits_{i=1}^n\overline\tau_{i}(\delta)}{|J_n|}=O(\delta^\alpha)\;\mbox{при}\;\delta\to 0,\quad 
\Big(\frac{|J_n|}{\sqrt{I_n}}\Big)^{1/(1+\alpha)}|\theta_n^*-\theta|\stackrel{p}{\to}0,\quad 0< \alpha\leq 1.
\end{eqnarray*}
В частности, в двух крайних случаях получаем следующие  ограничения на точность $\theta_n^*$:
$ \frac{|J_{n}|}{\sqrt{I_{n}}}|\theta_n^*-\theta|=O_p(1)$ при $\alpha=0$ и $ \frac{|J_n|^{1/2}}{I_n^{1/4}}|\theta_n^*-\theta|\stackrel{p}{\to}0 $ при $\alpha=1 $. Таким образом, если   $|J_n|/\sqrt{I_n}\sim \sqrt{n}$, то указанные ограничения на точность $\theta_n^*$ есть  либо  предположение о $\sqrt{n}$-ограниченности, либо об $n^{1/4}$-состоятельности  предварительной оценки. 
Отметим еще, что для  одношагового применения  изучаемых оценок
   при широких ограничениях  условие $ \frac{|J_n|^{1/2}}{I_n^{1/4}}|\theta_n^*-\theta|\stackrel{p}{\to}0 $ близко к минимальному достаточному условию на точность предварительной оценки
       (этот факт   доказан в  [\ref{2015-1?}]). 
Таким образом, асимптотическая нормальность одношаговых оценок установлена в достаточно широком  спектре ограничений на точность предварительной оценки $\theta_n^*$.
\hfill$\square$
\end{rem}%

\begin{rem}\label{p1-usl2}
Отметим, что  полученные в работе  условия асимптотической нормальности одношаговых   $M$-оценок не гарантируют не только состоятельность
 соответствующей  $M$-оценки,
 %
 но и даже ее  существование.
 Дело в том, что при некоторых
локальных и весьма слабых ограничениях (не гарантирующих даже существование $M$-оценки) в окрестности  $\theta$
с вероятностью, стремящейся к~$1$, существует и  единственен корень $\widetilde\theta_n(\theta)$  уравнения (\ref{p0-0+s}) --- ближайший к $\theta$ корень этого уравнения, обладающий требуемой точностью.
  Указанную
  случайную величину $\widetilde\theta_n(\theta)$, зависящую от $\theta$ и не являющуюся статистикой,   и приближают одношаговые $M$-оценки
 (подробнее объяснения отмеченному эффекту, а также библиографические ссылки,
 приведены  в [\ref{2015-2?}]). Аналогичное замечание имеет место и для  одношаговых взвешенных  $M$-оценок.

\end{rem}

\section{Приложения к задачам регрессии}
%

Рассмотрим некоторые статистические приложения полученных здесь результатов.  
Отметим, что различные примеры статистических  задач, в которых возникает $M$-оценивание в случае разнораспределенной выборки в схеме серий (и, как следствие, зачастую возникает и необходимость в применении одношаговых процедур) можно найти, например, в монографиях [\ref{2007-B}], [\ref{1997-H}], [\ref{1984-H}], [\ref{2006-J}], [\ref{2012-J1}]. 
%
    Здесь мы приведем   несколько примеров, относящихся к оцениванию в классических   задачах как  гомоскедастической, так и гетероскедастической  нелинейной регрессии.

{\bf 3.1.}
Пусть наблюдения $X_i$ 
  представимы в виде
 \begin{equation}  \label{p1-220}
 X_i=f_i(\theta)+\varepsilon_i\quad\mbox{при} \quad f_i(\cdot)=f(\cdot\;, a_{i}),\qquad i=1,\ldots,n,
 \end{equation}
 где  $f(\cdot,\cdot)$
 --- известная функция, при всех $i$ значение (скалярного или векторного) аргумента $a_{i}$  известно,
   погрешности $\{\varepsilon_i\}$ --- последовательность независимых
 случайных величин с нулевыми средними, а условие гомоскедастичности означает, что ${\bf D}\varepsilon_i=\sigma^2$ при всех $i$ (в противном случае говорят о гетероскедастичности модели). 

Один из основных 
методов нахождения оценки параметра $\theta$ в этой задаче --- метод наименьших квадратов (МНК) ---  состоит в минимизации суммы квадратов ошибок $\{\varepsilon_i\},$ т.е. МНК-оценка есть  
$\widetilde\theta_n=\arg\min\limits_\theta \sum\limits_{i=1}^n\big(X_i-f_i(\theta)\big)^2.$
При некоторых  условиях регулярности
отыскание  МНК-оценки сводится к решению уравнения 
\begin{equation}  \label{p1-226}
\sum\limits_{i=1}^nf'_i(t) (X_i-f_i(t))=0
 \end{equation}
Если  функция  $f_i(\theta)$ линейна по~$\theta$, 
 то МНК-оценка 
 считается в явном виде и  обладает рядом оптимальных свойств: является наилучшей линейной несмещенной оценкой, а в случае нормально распределенных ошибок ---  наилучшей и в классе всех несмещенных оценок (см., например,
[\ref{1981-Dem}], [\ref{1987-D1-2}], [\ref{1980-Seb}]). 

  Если же зависимость наблюдений от параметра нелинейная, то уравнение (\ref{p1-226}), как правило, не может быть решено аналитически.
   Более того, свойства линейных МНК-оценок не переносятся автоматически на нелинейный случай:  в общем случае, нелинейные МНК-оценки 
могут не существовать или не быть состоятельными 
  (см., например, [\ref{1981-Dem}]). Тем не менее, при выполнении некоторых условий регулярности,   нелинейные МНК-оценки 
 все же обладают некоторыми оптимальными свойствами по-крайней мере в асимптотике. Обоснованием точности нелинейного МНК в случае, когда ошибки имеют нормальное распределение, может служить совпадение в этом случае МНК-оценок  и оценок метода максимального правдоподобия. Подробную библиографию публикаций, связанных с изучением свойств нелинейных МНК-оценок, можно найти, например, в~[\ref{2003-S}]. Отметим, что при некоторых условиях регулярности  асимптотическая дисперсия МНК-оценки есть
 \begin{equation}\label{p1-223}
 {\sigma^2}\Big/{\sum\limits_{i=1}^n \big(f_i'(\theta)\big)^2}.
 \end{equation}

Таким образом, в нелинейном МНК  главная проблема
 касается вычислительной стороны: как отыскать оценку   $\widetilde\theta_n$, определяемую уравнением   (\ref{p1-226}), и являющуюся точкой глобального  минимума соответствующей функции.  Существует ряд приближенных численных методов решения этой задачи 
  и  выбор того или иного метода зависит от многих факторов, в том числе от вида функции $f_i(\theta)$.
 Зачастую проблема поиска нелинейной МНК-оценки связана с наличием  большого числа локальных минимумов у минимизируемой  функции, а потому
 при неудачном выборе начального приближения параметра итерационные процедуры  обнаруживают лишь локальный минимум, ближайший к этой стартовой точке.
 Указанный факт существенно затрудняет задачу  глобальной минимизации 
   (подробнее о вопросах приближенного поиска нелинейной МНК-оценки и возникающих при этом трудностях  см., например, 
   [\ref{1981-Dem}], [\ref{1989-D}],  [\ref{1987-D1-2}], [\ref{2003-S}]).

Наличие некоторой достаточно хорошей предварительной оценки $\theta_n^*$ позволяет,    используя  одношаговые процедуры, отыскать  в явном виде асимптотически  оптимальную оценку параметра $\theta$ и  избежать  указанные выше
вычислительные проблемы поиска нелинейной МНК-оценки. 

На примере  этой  задачи проиллюстрируем  еще раз  идеи, изложенные во введении.
С одной стороны, одношаговая $M$-оценка $\theta_{n}^{**},$ определенная в (\ref{p0-0s}) и приближающая состоятельное  решение $\widetilde\theta_n$ уравнения   (\ref{p1-226}), примет следующий вид:
 \begin{equation}  \label{p1-222}
 \theta_{n}^{**}=\theta_n^*+ \sum\limits_{i=1}^n
 (X_i-f_i(\theta_n^*))f'_i(\theta_n^*)\Big/\sum\limits_{i=1}^n\big[(f'_i(\theta_n^*))^2-(X_i-f_i(\theta_n^*))f''_i(\theta_n^*)\big].
 \end{equation}
 При выполнении некоторых условий  эта одношаговая $M$-оценка
  оказывается асимптотически нормальной с  асимптотической дисперсией из (\ref{p1-223}) (см.[\ref{2015-2?}] или следствие \ref{p1-cort1}).

 С другой стороны, в силу специфики уравнения (\ref{p1-226}),   с целью упростить одношаговую оценку и ослабить  условия на гладкость $f_i(\theta)$, вместо поиска решений и  их приближений для уравнения  (\ref{p1-226}), разумно рассматривать взвешенные $M$-оценки, определяемые   уравнением
\begin{equation}  \label{p1-2266}
\sum\limits_{i=1}^n h_i(\theta_n^*)M_i\big(t,X_i\big)=0 \quad \mbox{при}\quad h_i(t)=f'_i(t)\;\mbox{и}\; M_i(t,X_i)=X_i-f_i(t).
 \end{equation}
 В этом случае при некоторых дополнительных ограничениях  одношаговая взвешенная   $M$-оценка
 \begin{equation}  \label{p1-6222}
 \theta_{n,M,h}^{**}=\theta_n^*+ \sum\limits_{i=1}^n
 (X_i-f_i(\theta_n^*))f'_i(\theta_n^*)\Big/\sum\limits_{i=1}^n(f'_i(\theta_n^*))^2
 \end{equation}
  оказывается асимптотически нормальной с асимптотической дисперсией из  (\ref{p1-223})  (подчеркнем, что для функций $M_i(t,X_i)$ из (\ref{p1-2266}) функции  $h_i(t)$ из (\ref{p1-2266}) есть оптимальные функции $h_i^o(t)$, определенные в  предложении \ref{p1-t4-5}).
   Если при исследовании свойств  оценки $\theta_{n,M,h}^{**}$ из   (\ref{p1-6222}) воспользоваться
    теоремой \ref{p1-t1-2+}, то относительно гладкости функций $f_i(t)$ достаточно требовать, чтобы  функции $f_i'(t)$  удовлетворяли  лишь условию Гельдера.
Понятно, что при исследовании   оценки $\theta_{n}^{**}$ 
 требования на гладкость функций $f_i(t)$ оказываются сильнее  (см. определение (\ref{p1-222}) и следствие \ref{p1-cort1}).
   Таким образом,  одношаговые оценки   $\theta_{n}^{**}$ и $\theta_{n,M,h}^{**}$, определенные в (\ref{p1-222}) и  (\ref{p1-6222}),  асимптотически имеют    одну и ту же точность, при этом одношаговая взвешенная   $M$-оценка проще и является асимптотически нормальной при меньших условиях на гладкость
    функций $f_i(t)$.

  Отметим еще, что  если наблюдения $\{X_i\}$ имеют нормальные распределения, то информация Фишера $I_i(\theta)$, соответствующая наблюдению $X_i$, есть $I_i(\theta)=\big(f'_i(\theta)\big)^2/\sigma^2$; т.е. в этом случае асимптотическая дисперсия нелинейной МНК-оценки и одношаговых оценок
$\theta_{n}^{**}$ и $\theta_{n,M,h}^{**}$, определенных в  (\ref{p1-222}) и (\ref{p1-2266}),  обратно пропорциональна суммарной информации Фишера,  
соответствующей  $\{X_i\}:$
\begin{equation}  \label{p1-224}
{\sigma^2}\Big/{\sum\limits_{i=1}^n \big(f'_i(\theta)\big)^2}\equiv {1}\Big/{\sum\limits_{i=1}^n I_i(\theta)}.
 \end{equation}
Принимая во внимание обобщение неравенства Рао-Крамера на случай неоднородных наблюдений (см., например, [\ref{2007-B}]),
можно  сделать вывод  о точности соответствующих статистик.

\begin{rem} \label{p1-za}
Нетрудно видеть, что  сказанное выше очевидным образом переносится и на простейший  случай гетероскедастической регрессии, когда
\begin{equation}\label{p1-d2}
\forall i\quad {\bf E}\varepsilon_i=0,\quad{\bf D}\varepsilon_i=\sigma^2/w_i,
\quad \mbox{где\; числа}\; \{w_i\} - \mbox{известны}.
\end{equation}
Для этого  достаточно правую и левую часть уравнения (\ref{p1-220}) домножить на $\sqrt{w_i}$, что сводит  задачу к рассмотренной выше модели   гомоскедастической нелинейной регрессии. В частности, в приведенных выше рассуждениях нужно  вместо МНК-оценки рассматривать
 статистику
$\widetilde\theta_n=\arg\min\limits_\theta \sum\limits_{i=1}^n w_i\big(X_i-f_i(\theta)\big)^2$
(часто называемую взвешенной МНК-оценкой). В предположении (\ref{p1-d2}) для  одношаговой взвешенной $M$-оценки $\theta_{n,M,h}^{**}$ с оптимальной асимптотической дисперсией
  $I_{n,h^o}/J_{n,h^o}^2$ справедливы представления:
\begin{gather}\label{p1-222+}
\begin{split}
 \theta_{n,M,h}^{**}=\theta_n^*+\frac{ \sum\limits_{i=1}^n
 w_i(X_i-f_i(\theta_n^*))f'_i(\theta_n^*)}{\sum\limits_{i=1}^n w_i\big (f'_i(\theta_n^*)\big)^2},\quad
 I_{n,h^o}/J_{n,h^o}^2={\sigma^2}\Big/{\sum\limits_{i=1}^n w_i\big(f'_i(\theta)\big)^2}.
\end{split}
\end{gather}
\end{rem}

\begin{rem}\label{p1-zaz}
Предположим, что выполнено (\ref{p1-220}) и
$$
\forall i\quad {\bf E}\varepsilon_i=0,\quad{\bf E}\varepsilon_i=\sigma^2_i,
\quad \mbox{где\; значения }\; \{\sigma_i\} - \mbox{не \; известны}.
$$
В этом случае при $M_i(\theta,X_i)=X_i-f_i(\theta)$ оптимальные весовые функции $h_i^o(\theta)$, определяемые с точностью до константы равенством
$$h_i^o(\theta)={\bf E}M_i(\theta,X_i)/{\bf E}M_i^2(\theta,X_i)$$ (см. предложение \ref{p1-t4-5}), есть $h_i^o(\theta)=f'_i(\theta)/\sigma_i^2,$ т.е. оптимальные функции $h_i^o(\theta)$ не известны.
\end{rem}

Регрессионные функции $f(\cdot,\cdot)$ в  примерах \ref{p1-e2} и \ref{p1-en}
  взяты из монографии [\ref{1989-D}] и работ [\ref{2013-E}], [\ref{2014-S}].
 Подчеркнем, что
 в [\ref{1989-D}] эти  модели  приведены в качестве примеров, когда вычисление  МНК-оценки  затруднено по тем или иными причинам. 

\begin{ex}\label{p1-e2}
{\em
Пусть регрессионная модель задается соотношениями:
\begin{equation*}  \label{p1-233}
X_i=\sqrt{1+a_i\theta}+\varepsilon_i, \quad {\bf E}\varepsilon_i=0,\quad {\bf D}\varepsilon_i=\sigma^2/w_i,\quad i=1,\ldots,n,
\end{equation*}
где числа $\{a_i>0\}$ и $\{w_i>0\}$ известны. В   [\ref{2014-S}] в качестве альтернативы для  оптимальной  МНК-оценки   предложено (в случае $w_i\equiv 1$)  для параметра $\theta>0$
использовать статистику
\begin{equation}  \label{p1-234}
\theta_n^*=\sum\limits_{i=1}^n c_{i}X_i^2\Big/\sum\limits_{i=1}^n c_{i}a_i,\quad \mbox{где}\quad \sum\limits_{i=1}^n c_{i}=0.
\end{equation}
В [\ref{2014-S}] доказана  асимптотическая нормальность $\theta_n^*$ с асимптотической дисперсией $ \sum\limits_{i=1}^n c_{i}^2{\bf D}X_i^2\Big/\Big( \sum\limits_{i=1}^n c_{i}a_i\Big)^2$ и установлено, что   ни при каком  выборе констант $\{c_i\}$ 
эта величина не совпадает с асимптотической дисперсией МНК-оценки.

Используя  одношаговую оценку (\ref{p1-6222}),  точность оценивания из [\ref{2014-S}] может быть улучшена  (см. формулы (\ref{p1-235-}) и (\ref{p1-6222!}) при $w_i\equiv 1$).
Более того, оценку (\ref{p1-234}) можно  обобщить  на гетероскедастический случай (\ref{p1-d2}), полагая
\begin{equation}  \label{p1-234+}
\theta_n^*=\sum\limits_{i=1}^n c_{i}w_i(X_i^2-1)\Big/\sum\limits_{i=1}^n c_{i}w_ia_i,\quad \mbox{где}\quad \sum\limits_{i=1}^n c_{i}=0, \quad \sum\limits_{i=1}^n c_{i}w_ia_i\neq 0.
\end{equation}
Чтобы при соответствующих ограничениях доказать асимптотическую нормальность $\theta_n^*$ из (\ref{p1-234+}),
нужно с очевидными изменениями повторить доказательство аналогичного утверждения из  [\ref{2014-S}]. Подчеркнем, что оценка   (\ref{p1-234+}), как и  (\ref{p1-234}),  асимптотически нормальна, но не оптимальна.
Асимптотически оптимальная оценка из (\ref{p1-222+}) с асимптотической дисперсией взвешенной МНК-оценки
\begin{gather}   \label{p1-235-}
I_{n,h^o}/J_{n,h^o}^2=\sigma^2\big/\sum\limits_{i=1}^nw_i\big(f_i'(\theta)\big)^2\equiv 4\sigma^2\Big/\sum\limits_{i=1}^nw_ia_i^2(1+\theta a_i)^{-1}
 \end{gather}
 примет здесь следующий вид:
\begin{gather}\label{p1-6222!}
\theta_{n,M,h}^{**}=\theta_n^*+{2\sum\limits_{i=1}^n\big(X_i/\sqrt{1+a_i\theta_n^*}-1\big)w_ia_i}\Big/{\sum\limits_{i=1}^n
 w_ia_i^2(1+a_i\theta_n^*)^{-1}}.
 \end{gather}
 В качестве предварительной оценки $\theta_n^*$ можно использовать оценку   из~(\ref{p1-234+}).
~\hfill$\square$
}
\end{ex}

{\bf 3.2.}
С точки зрения приложений важную роль играют  гетероскедастические  регрессионные модели,  в которых дисперсии наблюдений зависят не только от номера наблюдения, но и от основного неизвестного параметра. В этом разделе считаем, что параметр $\theta$ нужно оценить по наблюдениям $X_1,\ldots,X_n,$ имеющим структуру (\ref{p1-220}), при этом на последовательность погрешностей $\{\varepsilon_i\}$
накладывается ограничение
\begin{equation}\label{p1-601}
{\bf E}\varepsilon_i=0,\qquad {\bf D}\varepsilon_i=\sigma^2/w_i(\theta),\qquad i=1,\ldots,n,
\end{equation}
где $\{w_i(\cdot)\}$ --- известные функции.

На первый взгляд, по аналогии с замечанием \ref{p1-za}, в рассматриваемой  ситуации было бы естественно также использовать оценки взвешенного МНК. Тем не менее,   метод наименьших квадратов  с весовыми
коэффициентами  не обеспечивает здесь нужной точности. Дело в том, что оптимальные веса метода наименьших квадратов должны в данной ситуации зависеть от параметра,  но если веса, определяющие оценку взвешенного МНК, зависят от $\theta,$ то такая оценка  в общем случае оказывается несостоятельной (см., например, [\ref{1971-F}], [\ref{1997-H}]).

Один из подходов оценивания параметра $\theta$ в данной ситуации состоит в использовании   $M$-оценки $\widetilde \theta_n$, определяемой  так называемым  уравнением {\it квази-правдопо\-добия}  (см., например, [\ref{1997-H}]): 
\begin{equation}\label{p1-602}
\sum\limits_{i=1}^nw_i(t)f_i'(t)(X_i-f_i(t))=0.
\end{equation}
 $M$-оценка $\widetilde \theta_n$ 
является наилучшей оценкой в классе статистик, являющихся решением (по $t$) уравнений вида
$\sum\limits_{i=1}^n h_i(t)(X_i-f_i(t))=0$ (см., например, [\ref{1997-H}], а также  предложение \ref{p1-t4-5}). 
Стоит отметить, что  $M$-оценка $\widetilde \theta_n$, определяемая уравнением  (\ref{p1-602}), 
 не является МНК-оценкой  (в случае, когда функции $w_i(t)$ не постоянны).

 Поиск оценки квази-правдоподобия 
 связан со всеми теми вычислительными проблемами, которые затрагивались  выше в разделе~3.1.
Взвешенную $M$-оценку $\widehat \theta_n,$ асимптотически имеющую точность оценки квази-правдоподобия, естественно определить как решение  уравнения
\begin{equation}  \label{p1-2260+}
\sum\limits_{i=1}^nh_i(\theta_n^*) M_i(t,X_i)=0\qquad\mbox{при }\quad h_i(t)=w_i(t)f_i'(t)\quad\mbox{и}\quad  M_i(t,X_i)=X_i-f_i(t).
 \end{equation}
  Одношаговая взвешенная  $M$-оценка $\theta_{n,M,h}^{**}$ из (\ref{p1-7+}) примет следующий вид:
\begin{equation}  \label{p1-2222}
 \theta_{n,M,h}^{**}=\theta_n^*+ {\sum\limits_{i=1}^n
 w_i(\theta_n^*)f'_i(\theta_n^*)(X_i-f_i(\theta_n^*))}\Big/{\sum\limits_{i=1}^nw_i(\theta_n^*)\big(f'_i(\theta_n^*)\big)^2}.
 \end{equation}
При этом в силу предложения \ref{p1-t4-5}  функции $h_i(\cdot)$ из (\ref{p1-2260+}) есть оптимальные  функции $h_i^o(\cdot)$ для $M_i(\cdot,X_i)$ из (\ref{p1-2260+}), а асимптотическая дисперсия $I_{n,h^o}/J_{n,h^o}^2$ оценки квази-правдоподобия (\ref{p1-602}) и оценок $\widehat\theta_n$,   $\theta_{n,M,h}^{**}$ из  (\ref{p1-2260+}) и (\ref{p1-2222}), определяется  соотношением
\begin{equation}  \label{p1-2223}
 I_{n,h^o}/J_{n,h^o}^2\equiv \sigma^2\Big/\sum\limits_{i=1}^n w_i(\theta)[f_i'(\theta)]^2.
 \end{equation}
Тем самым,  достигнута  та же точность (см. (\ref{p1-222+})), как и в случае, когда веса $w_i=w_i(\theta)$ были известны и не зависели от $\theta$ (см. также замечание \ref{p1-prim+}).

\begin{rem}\label{p1-prim}
Стоит отметить, что  взвешенная  $M$-оценка, имеющая асимптотическую дисперсию  (\ref{p1-2223}),
 может быть   определена и как решение следующего обобщения уравнения  (\ref{p1-2260+}):
\begin{eqnarray}  \label{p1-2260++}
\begin{split}
\sum\limits_{i=1}^nh_i(\theta_n^*) M_i(t,X_i)=0,\quad\mbox{где}\quad
h_i(t)=w_i(t)f_i'(t)g_i^{-1}(t),\;  M_i(t,X_i)=g_i(t)(X_i-f_i(t))
 \end{split}
 \end{eqnarray}
для некоторых функций $g_i(\cdot)$ (разумеется, $h_i(t)=h_i^o(t)$  для функций $M_i(t,X_i)$ из (\ref{p1-2260++})).  Удачный в том  или ином смысле  выбор функций  $g_i(t)$ зависит от специфики функций  $w_i(t)$ и $f_i(t)$, а уравнение (\ref{p1-2260++}) иногда может быть предпочтительнее, нежели уравнение  (\ref{p1-2260+})  (см. пример \ref{p1-e1}).  Далее, одношаговые  взвешенные $M$-оценки
можно строить исходя уже из уравнения  (\ref{p1-2260+}). Таким образом, за счет выбора  функций $g_i(\cdot)$  можно предложить  множество различных  явных оценок, эквивалентных в смысле асимптотической точности.
Кроме того, в [\ref{2015-2?}] предложены некоторые другие одношаговые оценки (так называемые оценки  скоринга), которые могут также определяться не однозначно и которые также эквивалентны предложенным здесь одношаговым взвешенным $M$-оценкам  в смысле асимптотической точности.
\hfill$\square$
\end{rem}

 \begin{rem}\label{p1-prim+}
   Метод квази-правдоподобия   не требует полной характеризации распределения: достаточно лишь знание первых двух моментов  наблюдений $\{X_i\}$ (как функций от неизвестных параметров).
В  случае же {\it нормально распределенных} наблюдений 
одношаговые приближения  при выполнении (\ref{p1-601}) разумно   находить  для оценки максимального правдоподобия, которая оказывается несколько точнее, нежели оценка квази-правдоподобия (см. 
[\ref{2015-2?}]).~\hfill$\square$
 \end{rem}


\begin{ex}\label{p1-en}
{\em
Пусть наблюдения $X_i$ имеют следующую структуру:
\begin{equation}  \label{p1-230}
X_i=a_i\theta+b_ig(\theta)+\varepsilon_i, \quad {\bf E}\varepsilon_i=0,\quad {\bf D}\varepsilon_i=\sigma^2/w_i(\theta),\quad i=1,\ldots,n,
\end{equation}
где $\{a_i\}$ и $\{b_i\}$ --- известные числовые последовательности.
 В
 [\ref{2013-E}] в качестве альтернативы трудно вычислимой   МНК-оценки 
 построена статистика
$\theta_n^*=\sum\limits_{i=1}^n c_{i}X_i\Big/\sum\limits_{i=1}^n c_{i}a_i,
$
где числа  $\{c_{i}\}$  таковы, что $\sum\limits_{i=1}^n c_{i}b_i=0$ и $\sum\limits_{i=1}^n c_{i}a_i\neq0.$ В
 [\ref{2013-E}]  доказано, что при некоторых дополнительных ограничениях оценка $\theta_n^*$ является асимптотически нормальной с асимптотической дисперсией $ \sigma^2\sum\limits_{i=1}^n c_{i}^2
 w_i^{-1}(\theta)
 \Big/\Big( \sum\limits_{i=1}^n c_{i}a_i\Big)^2$.

Используя обсуждаемый здесь подход  точность оценивания из [\ref{2013-E}]  можно улучшить,    при этом  статистика $\theta_n^*$ из [\ref{2013-E}]  может быть использована  в качестве предварительной  оценки для построения  одношаговых оценок.
 В частности,  одношаговая взвешенная $M$-оценка
$\theta_{n,M,h}$  из (\ref{p1-2222})  определяется равенством
\begin{gather*}   \label{p1-232}
\theta_{n,M,h}^{**}=\theta_n^*+\frac{\sum\limits_{i=1}^n w_i(\theta_n^*)\big(X_i-a_i\theta_n^*-b_ig(\theta_n^*)\big)\big(a_i+b_ig'(\theta_n^*)\big)}{
\sum\limits_{i=1}^nw_i(\theta_n^*)\big(a_i+b_ig'(\theta_n^*)\big)^2},
 \end{gather*}
  а асимптотическая дисперсия этой оценки, совпадающая с асимптотической дисперсией в некотором смысле оптимальной  оценки квази-правдоподобия,   есть $I_{n,h^o}/J_{n,h^o}^2\equiv\sigma^2\big/\sum\limits_{i=1}^nw_i(\theta)\big(a_i+b_ig'(\theta)\big)^2.
   $ Нетрудно видеть, что ни при каком выборе констант $\{c_{i}\}$ 
   для оценки $\theta_n^*$ из [\ref{2013-E}] не может быть достигнута эта точность.
%
%
\hfill$\square$
}
\end{ex}
\begin{ex} \label{p1-e1}
{\em
Рассмотрим одномерную гетероскедастическую модель Михаэлиса--Ментен, играющую важную роль в биохимии (см., например, 
  [\ref{2001-2}], [\ref{2011-M}]). В этом случае  наблюдения $\{X_i\}$
имеют следующую структуру:
\begin{equation}  \label{p1-227}
X_i={a_i}/{(1+b_i\theta)}+\varepsilon_i,\qquad{\bf E}\varepsilon_i=0,\quad {\bf D}\varepsilon_i=\sigma^2/w_i(\theta),\quad   i=1,\ldots,n,
 \end{equation}
где $\{a_i>0\}$ и $\{b_i>0\}$ --- известные числовые последовательности. 
Согласно (\ref{p1-2222}),
явную  асимптотически нормальную  оценку $\theta_{n,M,h}^{**}$ с минимальной асимптотической дисперсией
 $I_{n,h^o}/J_{n,h^o}^2$ можно определить соотношениями:
\begin{eqnarray}  \label{p1-2224}
\begin{split}
\theta_{n,M,h}^{**}=\theta_n^*-\frac{\sum\limits_{i=1}^n\Big(X_i-\frac{\displaystyle a_i}{\displaystyle 1+b_i\theta_n^*}\Big)\frac{\displaystyle w_i(\theta_n^*)a_ib_i}{\displaystyle (1+b_i\theta_n^*)^2}}{
\sum\limits_{i=1}^n\frac{\displaystyle w_i(\theta_n^*)a_i^2b_i^2}{\displaystyle (1+b_i\theta_n^*)^4}
},\quad
\frac{I_{n,h^o}}{J_{n,h^o}^2}= \frac{\sigma^2}{\sum\limits_{i=1}^n w_i(\theta)a_i^2b_i^2/(1+b_i\theta)^4}.
\end{split}
 \end{eqnarray}
В качестве оценки $\theta_n^*$ в (\ref{p1-2224}) можно использовать следующую статистику (асимптотически нормальную при широких ограничениях), построенную и исследованную в 
[\ref{2000-1}]:
\begin{equation}  \label{p1-227+}
\theta_n^*=\sum\limits_{i=1}^nc_i(a_i-X_i)\Big/\sum\limits_{i=1}^nc_ib_iX_i,
 \end{equation}
 где $\{c_i\}$ --- некоторые числа. Тот факт, что асимптотическая дисперсия оценки $\theta_n^*$ из (\ref{p1-227+}), равная
  $ \sum\limits_{i=1}^n c_i^2(1+b_i\theta)^2{\bf D}\varepsilon_i\Big/\Big( \sum\limits_{i=1}^n c_ia_ib_i(1+b_i\theta)^{-1}\Big)^2,$ при любом выборе констант $\{c_i\}$ больше, чем оптимальная асимптотическая дисперсия из (\ref{p1-2224}), установлен в [\ref{2000-1}].

Отметим еще, что в [\ref{2000-1}] 
 для модели  (\ref{p1-227}) 
 построена следующая  явная   асимптотически нормальная с асимптотической дисперсией  ${I_{n,h^o}}/{J_{n,h^o}^2}$ 
 из (\ref{p1-2224})  оценка: 
\begin{equation}  \label{p1-229}
\widetilde\theta_{n}^{**}=\frac{\sum\limits_{i=1}^n\frac{\displaystyle a_ib_iw_i(\theta_n^*)(a_i-X_i)}{\displaystyle (1+b_i\theta_n^*)^3}}{
\sum\limits_{i=1}^n\frac{\displaystyle a_ib_i^2w_i(\theta_n^*)X_i}{\displaystyle (1+b_i\theta_n^*)^3}}\equiv \theta_n^*-\frac{\sum\limits_{i=1}^n\Big(X_i-\frac{\displaystyle a_i}{\displaystyle 1+b_i\theta_n^*}\Big)\frac{\displaystyle w_i(\theta_n^*)a_ib_i}{\displaystyle (1+b_i\theta_n^*)^2}}{
\sum\limits_{i=1}^n\frac{\displaystyle w_i(\theta_n^*)a_ib_i^2X_i}{\displaystyle (1+b_i\theta_n^*)^3}
}.
\end{equation}
В силу приведенного тождества  статистика $\widetilde\theta_n^{**}$ отличается от одношаговой взвешенной $M$-оценки $\theta_{n,M,h}^{**}$ из (\ref{p1-2224}) знаменателем вычитаемой дроби.

Но как согласуется  статистика  $\widetilde\theta_n^{**}$ из (\ref{p1-229})
  с  методологией данной работы? Оказывается,  $\widetilde\theta_n^{**}$ является  одной из  оценок,  предлагаемых в данной статье.
  Статистика  $\widetilde\theta_n^{**}$   есть взвешенная $M$-оценка $\widehat\theta_n$, являющая решением уравнения (\ref{p1-2260++}) при $g_i(t)=1+b_it$: при так определенных функциях $g_i(t)$ уравнение    (\ref{p1-2260++}) явно  разрешается  относительно
$t$!   Более того,  $\widetilde\theta_n^{**}$ есть еще и  одношаговая взвешенная  $M$-оценка, построенная по указанному уравнению.
   Отметим еще, что в [\ref{2000-1}] при построении улучшенных  оценок не использовалась ни теория $M$-оценивания, ни  идея приближенного поиска $M$-оценок и их модификаций и метод Ньютона.   Использовались лишь некоторые эвристические соображения,  связанные со специальным  видом регрессионной  функции $f_i(\theta)$ из (\ref{p1-227}). Но эти соображения  не позволяют отыскать, например,  явные асимптотически оптимальные оценки в задачах  из примеров \ref{p1-e2} и~\ref{p1-en}, или найти явные оценки, эквивалентные оценке максимального правдоподобия. Предлагаемый
  здесь подход, в том числе  в замечаниях  \ref{p1-prim} и \ref{p1-prim+}, не ориентирован  на специальное представление регрессионной  функции 
   и, как следствие, более универсален. 
 \hfill$\square$

}
\end{ex}%
Достаточные условия  асимптотической нормальности одношаговых взвешенных  $M$-оценок,  указанных  в этом параграфе,
нетрудно получить из  теоремы   \ref{p1-t1-2+}, поэтому все точные условия мы здесь опускаем.

\section{Доказательства}
{\bf 4.1.}
Для вывода теорем \ref{p1-t1-2+} и \ref{p1-t21}  нам потребуется несколько вспомогательных обозначений и утверждений.

Положим $u_{n}^*=\theta_n^*-\theta$, \quad $\rho_{n}^M(u):=\displaystyle\sum\limits_{i=1}^nh_i(\theta)\big(M_i'(\theta+u,X_i)-M_i'(\theta,X_i)\big),$
\begin{eqnarray*}
\overline\rho_{n}^{M}(u):=\sum\limits_{i=1}^nh_i(\theta)\big(M_i(\theta+u,X_i)-M_i(\theta,X_i)-u M_i'(\theta)\big),
\\
\rho_{n}^{Mh}(u):=\sum\limits_{i=1}^n\big(h_i(\theta+u)-h_i(\theta)\big)\big(M_i'(\theta+u,X_i)-M_i'(\theta,X_i)\big),
\\
\overline\rho_{n}^{Mh}(u):=\sum\limits_{i=1}^n\big(h_i(\theta+u)-h_i(\theta)\big)\big(M_i(\theta+u,X_i)-M_i(\theta,X_i)-u M_i'(\theta)\big),
\\
\rho_{n}^{h}(u):=\displaystyle\sum\limits_{i=1}^n\big(h_i(\theta+u)-h_i(\theta)\big)M_i'(\theta,X_i).
\end{eqnarray*}
Нам также потребуются обозначения
$$\gamma_n(\delta):=\sum\limits_{i=1}^n|h_i(\theta)|\tau_i(\delta,X_i),\qquad
\gamma_n^h(\delta):=\delta^p\sum\limits_{i=1}^n\overline h_i\tau_i(\delta,X_i)$$ и соотношения
 ${\bf E}\gamma_n(\delta)=\sum\limits_{i=1}^n|h_i(\theta)|\overline\tau_i(\delta)$ и
 ${\bf E}\gamma_n^h(\delta)=\delta^p\sum\limits_{i=1}^n\overline h_i\overline\tau_i(\delta)$ при $\delta>0$.

\begin{lem} \label{p1-lem5}
Пусть $\alpha_n(\cdot): [0,\infty)\to [0,\infty)$ -- монотонно неубывающий случайный процесс, $\beta_n(\cdot)={\bf E}\alpha_n(\cdot)$ и
 $\beta_n(\eta_n)\stackrel{p}{\to}0$ для некоторой случайной величины  $\eta_n\geq 0$.
 Тогда  $\alpha_n(\eta_n)\stackrel{p}{\to}0$.
\end{lem}
Доказательство леммы \ref{p1-lem5}  приведено в  [\ref{2015-2?}].
\begin{lem}\label{p1-lem10j}
Пусть выполнено предположение $(R)$. Тогда
$$
\rho_n^{M}(u_n^*)/J_{n,h}\stackrel{p}{\to}0,  \qquad u_n^*\rho_{n}^M(u_n^*)/\sqrt{I_{n,h}}\stackrel{p}{\to}0,\qquad
\overline\rho_{n}^{M}(u_n^*)/\sqrt{I_{n,h}}
\stackrel{p}{\to}0.
$$
\end{lem}
Д~о~к~а~з~а~т~е~л~ь~с~т~в~о. Заметим прежде всего, что  при  выполнении $(R)$
 \begin{equation}\label{p1-119-s}
  |u_n^*|\sum\limits_{i=1}^n|h_i(\theta)|\overline\tau_i(|u_n^*|)/\sqrt{I_{n,h}}\stackrel{p}{\to }0\qquad\mbox{и}\qquad
\sum\limits_{i=1}^n|h_i(\theta)|\overline\tau_i(|u_n^*|)
 \big/J_{n,h} \stackrel{p}{\to }0.
 \end{equation}
Действительно, первая сходимость в (\ref{p1-119-s}) следует из  {\rm (\ref{p1-27!!})}. Кроме того, с учетом этой сходимости для  любого $\varepsilon>0$
 \begin{eqnarray}\label{p1-118s}
 {\bf P}\Big(\sum\limits_{i=1}^n\frac{|h_i(\theta)|\overline\tau_i(|u_n^*|)}{|J_{n,h}|}
 >\varepsilon\Big)= {\bf P}\Big(\sum\limits_{i=1}^n\frac{|h_i(\theta)|\overline\tau_i(|u_n^*|)}{|J_{n,h}|}>\varepsilon,\; \frac{|J_{n,h}|}{\sqrt{I_{n,h}}}|u_n^*|\geq 1\Big)+\nonumber\\ {\bf P}\Big(\sum\limits_{i=1}^n\frac{|h_i(\theta)|\overline\tau_i(|u_n^*|)}{|J_{n,h}|}>\varepsilon,\; \frac{|J_{n,h}|}{\sqrt{I_{n,h}}}|u_n^*|< 1\Big)\leq
  {\bf P}\Big(|u_n^*|\sum\limits_{i=1}^n\frac{|h_i(\theta)|\overline\tau_i(|u_n^*|)}{{I_{n,h}}}>\varepsilon\big)+\nonumber\\
   {\bf P}\Big(\sum\limits_{i=1}^n\frac{|h_i(\theta)|\overline\tau_i(|u_n^*|)}{|J_{n,h}|}>\varepsilon,\; |u_n^*|< \frac{\sqrt{I_{n,h}}}{|J_{n,h}|}\Big)\to 0.\qquad
 \end{eqnarray}
При выводе сходимости в (\ref{p1-118s}) мы также учли, что в силу монотонности функции $\sum\limits_{i=1}^n|h_i(\theta)\overline \tau_i(\cdot)/|J_{n,h}|$ и условий $(R_3)$ и $(R_4)$
начиная с некоторых $n$ в правой  части (\ref{p1-118s}) под знаком вероятности --- невозможное событие.
Далее, в силу $(R_2)$
$$
 M_i(\theta+\delta,X_i)-M_i(\theta,X_i)-\delta M_i'(\theta,X_i)=\delta\int\limits_0^1\big(M_i'(\theta+v\delta,X_i)-M_i'(\theta,X_i)\big)dv,
$$
а потому $\overline\rho_{n}^{M}(u)=u\displaystyle\int\limits_0^1\rho_{n}^{M}(uv)dv$ и справедливы оценки
\begin{equation}  \label{p1-120s}
\frac{|{\rho}_n^M(u)|}{|J_{n,h}|}\leq \frac{\gamma_n(|u|)}{|J_{n,h}|}, \quad \frac{|u||\rho_n^M(u)|}{\sqrt{I_{n,h}}} \leq \frac{|u|\gamma_n(|u|)}{\sqrt{I_{n,h}}},\quad \frac{|\overline\rho_n^M(u)|}{\sqrt{I_{n,h}}} \leq \frac{|u|\gamma_n(|u|)}{\sqrt{I_{n,h}}}.
\end{equation}
 Утверждения леммы следуют теперь из леммы \ref{p1-lem5} при $\eta_n=|u_n^*|$, если положить   $\alpha_n(|u|)=\gamma_n(|u|)/|J_{n,h}|$  или   $\alpha_n(|u|)=|u|\gamma_n(|u|)/\sqrt{I_{n,h}}$ и учесть сходимости (\ref{p1-119-s}).
\hfill$\square$.

\begin{lem}\label{p1-lem11}
Если выполнены условия $(R)$, то
$$
\frac{\rho_n^{Mh}(u_n^*)}{J_{n,h}}\stackrel{p}{\to}0,  \quad \frac{u_n^*\rho_{n}^{Mh}(u_n^*)}{\sqrt{I_{n,h}}}\stackrel{p}{\to}0,\quad
\frac{\overline\rho_{n}^{Mh}(u_n^*)}{\sqrt{I_{n,h}}}
\stackrel{p}{\to}0,\quad \frac{\rho_n^{h}(u_n^*)}{J_{n,h}}\stackrel{p}{\to}0.
$$
  \end{lem}
Д~о~к~а~з~а~т~е~л~ь~с~т~в~о первых трех соотношений  аналогично выводу утверждений леммы \ref{p1-lem10j}, поэтому подробные рассуждения мы опускаем.
 В силу $(R_2)$
\begin{equation}\label{p1-53w}
|\rho_n^h(u)|/|J_{n,h}|\leq |u|^p\sum\limits_{i=1}^n\overline h_i|M_i'(\theta,X_i)|/|J_{n,h}|,
\end{equation}
 тем самым последнее утверждение леммы  следует из (\ref{p1-j27}) и леммы   \ref{p1-lem5} при $\eta_n=|u_n^*|$, если  в качестве $\alpha_n(|u|)$ рассматривать правую часть в (\ref{p1-53w}).
\hfill$\square$

%
%
%

 Д~о~к~а~з~а~т~е~л~ь~с~т~в~о теоремы \ref{p1-t1-2+}.
  Поскольку в силу  (\ref{p1-119-s}) c вероятностью, стремящейся к  $1$  величина $
  \sum\limits_{i=1}^n|h_i(\theta)|\overline\tau_i(|\theta_n^*-\theta|)
 \big/|J_{n,h}|$ конечна, то ${\bf P}(\theta_n^*\in \Theta)\to 1.$
 Следовательно, согласно условию $(R_2)$, с вероятностью, стремящейся к~$1$ определены  функции  $\{h_i(\theta_n^*)\}$, $\{M_i(\theta_n^*)\}$ и $\{M_i'(\theta_n^*)\}$, участвующие в представлении (\ref{p1-7+})  для оценки $\theta_{n,M,h}^{**}$.
В силу обозначений в начале раздела 4.1, справедливо равенство
 $\sum\limits_{i=1}^nh_i(\theta_n^*)M'_i(\theta_n^*,X_i)=\sum\limits_{i=1}^nh_i(\theta)M'_i(\theta,X_i)
 +\rho_n^M(u_n^*)+\rho_n^{Mh}(u_n^*)+\rho_n^h(u_n^*),$
 а потому, с учетом  лемм  \ref{p1-lem10j}, \ref{p1-lem11} и условия  $(R_3)$, выполнено
  ${\bf P}\Big(\sum\limits_{i=1}^nh_i(\theta_n^*)M'_i(\theta_n^*,X_i)\neq 0\Big)\to 1.$
  Таким образом, согласно определению (\ref{p1-7+}), статистика  $\theta_{n,M,h}^{**}$ определена с вероятностью, стремящейся к~$1$.
Далее,  из определения (\ref{p1-7+})   
имеем
\begin{equation}                                                                       \label{p1-100+}
\theta_{n,M,h}^{**}-\theta=-\frac{\sum\limits_{i=1}^n h_i(\theta_n^*)M_i(\theta_n^*,X_i)-(\theta_n^*-\theta)\sum\limits_{i=1}^n h_i(\theta_n^*)M_i'(\theta_n^*,X_i)}{\sum\limits_{i=1}^n h_i(\theta_n^*)M_i'(\theta_n^*,X_i)}.
\end{equation}
Используя представление (\ref{p1-100+}) и
обозначения, введенные перед леммой \ref{p1-lem5},
получаем следующее тождество:$\quad \frac{\displaystyle J_{n,h}}{\displaystyle \sqrt{I_{n,h}}}(\theta_{n,M,h}^{**}-\theta)= $
\begin{eqnarray}                                                                     \label{p1-102+}
-\frac{\Big(\sum\limits_{i=1}^n h_i(\theta) M_i(\theta,X_i)-u_n^*\rho_n^M(u_n^*)
+\overline\rho_n^M(u_n^*)-u_n^*\rho_n^{Mh}(u_n^*)
+\overline\rho_n^{Mh}(u_n^*)\Big)/\sqrt{I_{n,h}}+\delta_n^h}
{\sum\limits_{i=1}^n h_i(\theta)M_{i}'(\theta,X_i)/J_{n,h}+\big(\rho_n^M(u_n^*)+\rho_n^{Mh}(u_n^*)+\rho_n^h(u_n^*)\big)/J_{n,h}}.
  \end{eqnarray}
    Чтобы теперь  из тождества (\ref{p1-102+}) извлечь асимптотическую нормальность (\ref{p1-24-}) оценки
 $\theta_{n,M,h}^{**}$, остается воспользоваться утверждениями лемм \ref{p1-lem10j}, \ref{p1-lem11} и
 условиями  (\ref{p1-5-}) и   (\ref{p1-j41}).\hfill$\square$

\begin{lem} \label{p1-lem2+}
Пусть 
выполнены
 условия теоремы {\rm\ref{p1-t21}}. Тогда
\begin{equation*}
\delta_{n}^M:=
\Big(\big(\sum\limits_{i=1}^nh_i^2(\theta_{n,M,h}^{**})M_i^2(\theta_{n,M,h}^{**},X_i)\big)^{1/2}-
\big(\sum\limits_{i=1}^nh_i^2(\theta)M_i^2(\theta,X_i)\big)^{1/2}
\Big)/\sqrt{I_{n,h}}\stackrel{p}{\to}0.
\end{equation*}
\end{lem}
Д~о~к~а~з~а~т~е~л~ь~с~т~в~о. 
 Положим 
  $u_n^{**}:=\theta_{n,M,h}^{**}-\theta$ и
условимся у функций
$M_i(\cdot,X_i)$.
 опускать второй аргумент.
Нам потребуются соотношения
 \begin{equation}  \label{p1-117-s}
J_{n,h}u_n^{**}/\sqrt{I_{n,h}}=O_p(1),\quad u_n^{**}\stackrel{p}{\to}0,\quad \gamma_{n}(|u_n^{**}|)/J_{n,h}
\stackrel{p}{\to}0,\quad \gamma_{n}^h(|u_n^{**}|)/J_{n,h}
\stackrel{p}{\to}0.
\end{equation}
Первое из этих соотношений  следует из теоремы \ref{p1-t1-2+}, второе --- из  условия $(R_3)$ и первого соотношения.  Остальные -- из леммы \ref{p1-lem5} при $\eta_n=|u_n^{**}|$ и
 $\alpha_n(\cdot)=\gamma_{n}(\cdot)/|J_{n,h}|$ или $\alpha_n(\cdot)=\gamma_{n}^h(\cdot)/|J_{n,h}|$, если только еще учесть
 условие $(R_4)$. Далее,  имеем
\begin{eqnarray}\label{p1-54w}
\sqrt{I_{n,h}}|\delta_{n}^M|\leq \Big(\sum\limits_{i=1}^n\big(h_i(\theta_{n,M,h}^{**})M_i(\theta_{n,M,h}^{**})-h_i(\theta)M_i(\theta)\big)^2\Big)^{1/2}=
\delta_{n1}^M+\delta_{n2}^M+\delta_{n3}^M,
\\
\delta_{n1}^M:=\Big(\sum\limits_{i=1}^nh_i^2(\theta)\big(M_i(\theta_{n,M,h}^{**})-M_i(\theta)\big)^2\Big)^{1/2}, \nonumber\\
\delta_{n2}^M:=\Big(\sum\limits_{i=1}^n\big(h_i(\theta_{n,M,h}^{**})-h_i(\theta)\big)^2\big(M_i(\theta_{n,M,h}^{**})-M_i(\theta)\big)^2\Big)^{1/2},\nonumber\\
\delta_{n3}^M:=\Big(\sum\limits_{i=1}^n\big(h_i(\theta_{n,M,h}^{**})-h_i(\theta)\big)^2M_i^2(\theta)\Big)^{1/2}.\nonumber
\end{eqnarray}
Оценим отдельно каждую из введенных величин. В силу предположения  $(R_2)$
\begin{eqnarray*}
\delta_{n1}^M\leq \Big(\sum\limits_{i=1}^nh_i^2(\theta)\Big(\int\limits_{\theta}^{\theta_{n,M,h}^{**}}M_i'(t)dt\Big)^2\Big)^{1/2}\leq
\sqrt{2}\Big(\sum\limits_{i=1}^nh_i^2(\theta)\Big(\int\limits_{\theta}^{\theta_{n,M,h}^{**}}\big(M_i'(t)-M_i'(\theta)\big)dt\Big)^2\Big)^{1/2}+\nonumber\\
\sqrt{2}\Big(\sum\limits_{i=1}^nh_i^2(\theta)\Big(\int\limits_{\theta}^{\theta_{n,M,h}^{**}}M_i'(\theta)dt\Big)^2\Big)^{1/2}\leq
\sqrt{2}|u_n^{**}|\Big(\sum\limits_{i=1}^nh_i^2(\theta)\Big(\sup\limits_{t:|t-\theta|\leq |u_n^{**}|}\big|M_i'(t)-M_i'(\theta)\big|\Big)^2\Big)^{1/2}+\nonumber\\
\sqrt{2}|u_n^{**}|\Big(\sum\limits_{i=1}^nh_i^2(\theta)\big(M_i'(\theta)\big)^2\Big)^{1/2}\leq
\sqrt{2}|u_n^{**}|\gamma_{n}(|u_n^{**}|)+\sqrt{2}|u_n^{**}|\Big(\sum\limits_{i=1}^nh_i^2(\theta)\big(M_i'(\theta)\big)^2\Big)^{1/2}.
\end{eqnarray*}
Следовательно, используя  соотношения (\ref{p1-117-s})  условие  (\ref{p1-8-}), получаем
 \begin{equation*}
 \frac{\delta_{n1}^M}{\sqrt{I_{n,h}}}\leq \sqrt{2}\frac{|J_{n,h}|}{\sqrt{I_{n,h}}}|u_n^{**}|\frac{\gamma_{n}(|u_n^{**}|)}{|J_{n,h}|}+\sqrt{2}
\frac{|J_{n,h}|}{\sqrt{I_{n,h}}}|u_n^{**}|\Big(\sum\limits_{i=1}^nh_i^2(\theta)\big[M_i'(\theta)\big]^2/J_{n,h}^2\Big)^{1/2}\stackrel{p}{\to}0.
\end{equation*}
Проводя аналогичные рассуждения, имеем
\begin{equation*}
\frac{\delta_{n2}^M}{\sqrt{I_{n,h}}}\leq \sqrt{2}\frac{|J_{n,h}|}{\sqrt{I_{n,h}}}|u_n^{**}|\frac{\gamma_{n}^h(|u_n^{**}|)}{|J_{n,h}|}+\sqrt{2}
\frac{|J_{n,h}|}{\sqrt{I_{n,h}}}|u_n^{**}|^{1+p}\Big(\sum\limits_{i=1}^n\overline h_i^2\big[M_i'(\theta)\big]^2/J_{n,h}^2\Big)^{1/2}\stackrel{p}{\to}0,
\end{equation*}
\begin{equation*}
\frac{\delta_{n3}^M}{\sqrt{I_{n,h}}}\leq \Big(\frac{|J_{n,h}||u_n^{**}|}{\sqrt{I_{n,h}}}\Big)^p\Big(\sum\limits_{i=1}^n\overline h_i^2\big[M_i'(\theta)\big]^2\big/(J_{n,h}^{2p}I_{n,p}^{1-p})\Big)^{1/2}\stackrel{p}{\to}0.
\end{equation*}
Три последние сходимости, вместе с оценкой (\ref{p1-54w}), завершают доказательство.
\hfill $\square$

  Д~о~к~а~з~а~т~е~л~ь~с~т~в~о теоремы \ref{p1-t21}.
Поскольку ${\bf P}(\theta_n^*\in \Theta)\to 1$ (см. доказательство теоремы \ref{p1-t1-2+}) и ${\bf P}(\theta_{n,M,h}^{**}\in \Theta)\to 1$ в силу
(\ref{p1-117-s}) и $(R_3)$,
 то  с вероятностью, стремящейся к $1$ определены  все величины в правой части равенства в (\ref{p1-9-s}). Кроме того,
  поскольку в силу   леммы   \ref{p1-lem2+}  и условия  (\ref{p1-8-})
  \begin{equation}\label{p1-51w}
  \frac{\Big(\sum\limits_{i=1}^n h_i^2(\theta_{n,M,h}^{**})M_i^2(\theta_{n,M,h}^{**},X_i)\Big)^{1/2}}{\sqrt{I_{n,h}}}\equiv \frac{\Big(\sum\limits_{i=1}^n h_i^2(\theta)M_i^2(\theta,X_i)\Big)^{1/2}}{\sqrt{I_{n,h}}}+\delta_{n}^M\stackrel{p}{\to}1,
  \end{equation}
  то  знаменатель в правой части равенства в (\ref{p1-9-s}) обращается в нуль с вероятностью, стремящейся к нулю.
 Таким образом, с учетом определения (\ref{p1-9-s}), статистика  $d_{n,M,h}^{*}$ определена с вероятностью, стремящейся к $1$.
Далее, учитывая представления  (\ref{p1-100+}),  (\ref{p1-102+}) и  обозначения
  из (\ref{p1-9-s}), имеем $d_{n,M,h}^*(\theta_{n,M,h}^{**}-\theta)=$
\begin{eqnarray*}
-\frac{\Big(\sum\limits_{i=1}^n h_i(\theta) M_i(\theta,X_i)-u_n^*\rho_n^M(u_n^*)
+\overline\rho_n^M(u_n^*)-u_n^*\rho_n^{Mh}(u_n^*)
+\overline\rho_n^{Mh}(u_n^*)\Big)/\sqrt{I_{n,h}}+\delta_n^h}
{\big(\sum\limits_{i=1}^n h_i^2(\theta)M_i^2(\theta,X_i)\big)^{1/2}/\sqrt{I_{n,h}}+\delta_{n}^M}.
\end{eqnarray*}
Тот факт, что  числитель в правой части этого представления  по распределению сходится к стандартному нормальному закону, установлен при выводе теоремы \ref{p1-t1-2+}. Отмеченная сходимость вместе с  (\ref{p1-51w})  доказывает (\ref{p1-9-ss}).
\hfill$\square$

{\bf 4.2.}
Д~о~к~а~з~а~т~е~л~ь~с~т~в~о предложения \ref{p1-pre2}. Положим
$$\theta_{ni}^*= {\bf E}_{\neq i}\theta_n^*,\qquad\theta_{nij}^*:={\bf E}_{\neq
 i}{\bf E}_{\neq j}\theta_n^*,\qquad  i,j=1,\ldots,n,\; j\neq i. $$
 В дальнейшем не единожды, больше этого не оговаривая, будет использован тот факт, что    величина $\theta_{ni}^*$ не
зависят от наблюдения $X_i$, а величина $\theta_{nij}^*$ не зависят еще и от  $X_j$.
Заметим прежде всего, что  при всех   $i$ и $j\neq i$ справедливы оценки
\begin{equation} \label{p1-f1}
{\bf E}\big|\theta_{ni}^*-\theta\big|^2\leq  4{\bf E}\big|\theta_n^*-\theta \big|^2
\qquad {\bf E}\big|\theta_{ni}^*-\theta_{nij}^*
 \big|^2\leq {\bf E}\big|\theta_{n}^*-\theta_{nj}^*
 \big|^2.
\end{equation}
Действительно,
\begin{equation} \label{p1-f2}
{\bf E}\big|\theta_{ni}^*-\theta\big|^2 \leq
2{\bf E}\big|\theta_{n}^*-\theta\big|^2+2{\bf E}\big|\theta_{n}^*-\theta_{ni}^*\big|^2.
\end{equation}
Кроме того, с учетом  формулы полной вероятности,
\begin{equation} \label{p1-f3}
{\bf E}\big|\theta_{n}^*-\theta_{ni}^*
\big|^2={\bf E}{\bf E}_{\neq
i}\big|\theta_n^*-{\bf E}_{\neq i}\theta_n^*\big|^2\leq {\bf E}{\bf E}_{\neq i}\big|\theta_n^*-h_{\neq i} \big|^2,
\end{equation}
где величина $h_{\neq i}$ не зависит от $X_i$. В частности,
полагая $ h_{\neq i}=\theta,$ из (\ref{p1-f2}) и (\ref{p1-f3}),   извлекаем первое утверждение в (\ref{p1-f1}).
Еще раз используя  формулу полной вероятности, получаем следующую цепочку соотношений:
\begin{equation*}
{\bf E}\big|\theta_{ni}^*-\theta_{nij}^*
 \big|^2={\bf E}\big|{\bf E}_{\neq i}\theta_n^*-{\bf E}_{\neq i}{\bf E}_{\neq j}
 \theta_n^* \big|^2\leq{\bf E}{\bf E}_{\neq i}\big|\theta_n^*-{\bf E}_{\neq j}
 \theta_n^* \big|^2={\bf E}\big|\theta_n^*- \theta_{nj}^*
 \big|^2,
\end{equation*}
 которая доказывает  второе утверждение в (\ref{p1-f1}).
%

Имеем
\begin{equation}  \label{p1-j45}
\delta_{n}^h=\frac{\sum\limits_{i=1}^n\big(h_i(\theta_n^*)-h_i(\theta_{ni}^*)\big)M_i(\theta,X_i)}{\sqrt{I_{n,h}}}+
\frac{\sum\limits_{i=1}^n\big(h_i(\theta_{ni}^*)-h_i(\theta)\big)M_i(\theta,X_i)}{\sqrt{I_{n,h}}}\equiv\delta_{n,1}^h+\delta_{n,2}^h.
\end{equation}
В условиях предложения  $
\sqrt{I_{n,h}}{\bf E}|\delta_{n,1}^h|\leq $
$$\sum\limits_{i=1}^n\overline h_i{\bf E}\big(|\theta_n^*-\theta_{ni}^*|^p|M_i(\theta,X_i)|\big)\leq
\sum\limits_{i=1}^n\overline h_i\big({\bf E}|\theta_n^*-\theta_{ni}^*|^2\big)^{p/2}\big({\bf E}|M_i(\theta,X_i)|^2\big)^{1/2}.
$$
 Следовательно, с учетом  условия (\ref{p1-j29}),   ${\bf E}|\delta_{n,1}^h|\to 0,$ а потому, в силу неравенства Чебышева с первым моментом,   $\delta_{n,1}^h\stackrel{p}{\to}0.$

Покажем, что $\delta_{n,2}^h\stackrel{p}{\to}0.$
Нам потребуются следующие оценки:
\begin{eqnarray}  \label{p1-j46}
{\bf E}\big((h_i(\theta_{ni}^*)-h_i(\theta))^2M_i^2(\theta,X_i)\big)=
{\bf E}\big(h_i(\theta_{ni}^*)-h_i(\theta)\big)^2{\bf E}M_i^2(\theta,X_i)
\leq\nonumber\\ \overline h_i^2{\bf E}\big(|\theta_{ni}^*-\theta|^2\big)^p{\bf E}M_i^2(\theta,X_i)\leq 4^p{\bf E}\big(|\theta_{n}^*-\theta|^2\big)^p\overline h_i^2{\bf E}M_i^2(\theta,X_i),
\end{eqnarray}
\begin{eqnarray}  \label{p1-j47}
{\bf E}\Big[\big(h_i(\theta_{ni}^*)-h_i(\theta)\big)\big(h_j(\theta_{nj}^*)-h_j(\theta)\big)M_i(\theta,X_i)M_j(\theta,X_j)\Big]=\nonumber\\
{\bf E}\Big[\big(h_i(\theta_{ni}^*)-h_i(\theta_{nij}^*)+h_i(\theta_{nij}^*)-h_i(\theta)\big)\big(h_j(\theta_{nj}^*)-h_j(\theta_{nji}^*)+
h_j(\theta_{nji}^*)-h_j(\theta)\big)\times \nonumber\\M_i(\theta,X_i)M_j(\theta,X_j)\Big]=
{\bf E}\Big[\big(h_i(\theta_{ni}^*)-h_i(\theta_{nij}^*)\big)\big(h_j(\theta_{nj}^*)-h_j(\theta_{nji}^*)\big)M_i(\theta,X_i)M_j(\theta,X_j)\Big]+
\nonumber\\
{\bf E}\Big[\big(h_i(\theta_{ni}^*)-h_i(\theta_{nij}^*)\big)\big(h_j(\theta_{nji}^*)-h_j(\theta)\big)M_i(\theta,X_i)M_j(\theta,X_j)\Big]+
\nonumber\\
{\bf E}\Big[\big(h_i(\theta_{nij}^*)-h_i(\theta)\big)\big(h_j(\theta_{nj}^*)-h_j(\theta_{nji}^*)\big)M_i(\theta,X_i)M_j(\theta,X_j)\Big]+
\nonumber\\
{\bf E}\Big[\big(h_i(\theta_{nij}^*)-h_i(\theta)\big)\big(h_j(\theta_{nji}^*)-h_j(\theta)\big)M_i(\theta,X_i)M_j(\theta,X_j)\Big]=\nonumber\\
{\bf E}\Big[\big(h_i(\theta_{ni}^*)-h_i(\theta_{nij}^*)\big)\big(h_j(\theta_{nj}^*)-h_j(\theta_{nji}^*)\big)M_i(\theta,X_i)M_j(\theta,X_j)\Big]\leq \nonumber\\
\left({\bf E}\big(h_i(\theta_{ni}^*)-h_i(\theta_{nij}^*)\big)^2{\bf E}M_i^2(\theta,X_i)\right)^{1/2}
\left({\bf E}\big(h_j(\theta_{nj}^*)-h_j(\theta_{nji}^*)\big)^2{\bf E}M_j^2(\theta,X_j)\right)^{1/2}\leq \nonumber\\
\overline h_i\overline h_j\big({\bf E}|\theta_{n}^*-\theta_{ni}^*|^2\big)^{p/2}\big({\bf E}|\theta_{n}^*-\theta_{nj}^*|^2\big)^{p/2}\big({\bf E}M_i^2(\theta,X_i)\big)^{1/2}\big({\bf E}M_j^2(\theta,X_j)\big)^{1/2}.\nonumber
\end{eqnarray}
При выводе этих соотношений были учтены еще  неравенства из  (\ref{p1-f1}). Следовательно,
 с учетом условий (\ref{p1-j28}) и (\ref{p1-j29})
\begin{eqnarray*}
{\bf E}|\delta_{n,2}^h|^2=\sum\limits_{i=1}^n {\bf E}\big((h_i(\theta_{ni}^*)-h_i(\theta))^2M_i^2(\theta,X_i)\big)\big/I_{n,h}+
\nonumber\\
\sum\limits_{i,j=1, i\neq j}^n {\bf E}\Big[\big(h_i(\theta_{ni}^*)-h_i(\theta)\big)\big(h_j(\theta_{nj}^*)-h_j(\theta)\big)M_i(\theta,X_i)M_j(\theta,X_j)\Big]\big/I_{n,h}\leq \nonumber\\
\leq 4^p{\bf E}\big(|\theta_{n}^*-\theta|^2\big)^p\sum\limits_{i=1}^n\overline h_i^2{\bf E}M_i^2(\theta,X_i)/I_{n,h}+\nonumber\\
\left(\sum\limits_{i=1}^n\overline h_i\big({\bf E}|\theta_{n}^*-\theta_{ni}^*|^2\big)^{p/2}\big({\bf E}M_i^2(\theta,X_i)\big)^{1/2}\big/\sqrt{I_{n,h}}\right)^{2}
\to 0
\end{eqnarray*}
и доказательство предложения завершено.
\hfill$\square$

%
%
%
%

Д~о~к~а~з~а~т~е~л~ь~с~т~в~о предложения \ref{p1-t4-5}. Соотношение (\ref{p1-391}) нетрудно извлечь из неравенства Коши-Буняковского.  Положим $K_n:=\sum\limits_{i=1}^n\big(h_i^{o}(\theta)\big)^2$. %
Поскольку по условиям  теоремы   при каждом $i$ величины
$h_i(\theta)$ и $h_i^{o}(\theta)$ имеют одинаковый
знак, то справедливы  неравенства (\ref{p1-394}). В частности, при
каждом $i$ значение $h_i(\theta)$ лежит между
$h h_i^{o}(\theta)$ и $Hh_i^{o}(\theta)$. Значит для
величин
$d_i:=(h_i(\theta)-h h_i^{o}(\theta))(h_i(\theta)-
Hh_i^{o}(\theta))$, $ i=1,\ldots,n,$ выполнено \ $d_i\le0$. Из
этого факта с учетом обозначений из 
 $(R_3)$ получаем, что
$I_{n,h}-(H+h)J_{n,h}+hHK_n=\sum\limits_{i=1}^n d_i{\bf E}M_i^2(\theta,X_i)\le0. $
А поскольку функция $f(x):=f(x)=x^{-1}-Hh(H+h)^{-1}K_nx^{-2}$ достигает своего максимального значения при
   $x=2Hh(H+h)^{-1}K_n$, то
\begin{eqnarray*}
{I_{n,h}}/{J_{n,h}^2} \le
{(H+h)J_{n,h}^{-1}-HhK_n}{J_{n,h}^{-2}}=(H+h)f(J_{n,h})\le\nonumber\\
(H+h)f\Big(\frac{2HhK_n}{H+h}\Big)=\frac{(H+h)^2}{4HhK_n}
=\frac{I_{n,h^o}}{J_{n,h^o}^2}\frac{(H+h)^2}{4Hh}.
\end{eqnarray*}
 Значит,
$\frac{\displaystyle {I_{n,h}}/{J_{n,h}^2}}{\displaystyle{I_{n,h^o}}/{J_{n,h^o}^2}}
\le \frac{\displaystyle H+h}{\displaystyle 2\sqrt{Hh}} =\frac{\displaystyle H/h+1}{\displaystyle 2\sqrt{H/h}}
=1+\frac{\displaystyle(\sqrt{H/h}-1)^2}{\displaystyle 2\sqrt{H/h}}\leq \sqrt{\frac{\displaystyle H}{\displaystyle h}}
$
и доказательство предложения завершено.
 \hfill$\square$

\vspace{15mm}
\centerline{\Large \textbf{Литература}}
\addcontentsline{toc}{section}{Литература}

\begin{enumerate}

\item\label{2007-B}
{\em Боровков А.А. } Математическая статистика. Москва: Физматлит,
 2007.

%

\item\label{1988-G} {\em Гуревич В.А.} О взвешенных М-оценках в нелинейной регрессии.
 --- Теория  вероятностей и ее применения, 1988, Т.~33, 2, 421-424.
%
%

\item\label{1981-Dem} {\em Демиденко Е.З.} Линейная и нелинейная регрессия.  М.: Финансы и статистика, 1981.

\item\label{1989-D} {\em Демиденко Е.З.} Оптимизация и регрессия. М.: Наука, 1989.

\item\label{1987-D1-2} {\em Дрейпер Н., Смит Г.} Прикладной регрессионный анализ. Кн. 1,2.  М.: Финансы и статистика. 1987.

%
%
%

\item\label{2013-E}  {\em  Ермоленко К.В., Саханенко А.И. } Явные асимптотически нормальные оценки неизвестного  параметра частично-линейной регрессии. ---
 Сиб. электрон. матем. изв., 2013, Т.10, 719-726.

\item\label{1975-Z}
{\em Закс Ш.} Теория статистических выводов. М.:Мир, 1975.
%

\item\label{1991-L}
{\em Леман Э.}  Теория точечного оценивания. М.: Наука, 1991.

\item\label{2015-1?}
{\em  Линке Ю.Ю.} Об уточнении одношаговых оценок Фишера в случае медленно сходящихся предварительных  оценок. --- Теор. вероятн. и ее примен., 2015, Т.60, вып. 1, 80-98.

\item\label{2015-2?}
{\em  Линке Ю.Ю.} Асимптотические свойства  одношаговых $M$-оценок для разнораспределенных  наблюдений.  (сдана в печать) ArXiv:  1503.03393v8

\item\label{2011-1}
{\em  Линке Ю.Ю.} Об асимптотике распределения двухшаговых
статистических оценок. ---  Сиб. Матем. Жур., 2011,
Т.~52. \No~4, 841-860.

%
%

\item\label{2013-2}
{\em  Линке Ю.Ю., Саханенко А.И.} Об асимптотике распределения одного
класса двухшаговых статистических оценок многомерного параметра. ---
Математические труды, 2013, Т.~16. \No~1, 89-120.

\item\label{2014-1}
    {\em Линке Ю.Ю., Саханенко А.И.} Об условиях асимптотической нормальности одношаговых оценок Фишера для однопараметрических семейств распределений. ---   Сиб. электрон. матем. изв., 2014, Т.11, 464-475.

\item\label{2000-1}
{\em Линке Ю.Ю. Саханенко А.И.}  Асимптотически нормальное
оценивание параметра в задаче дробно -- линейной регрессии.
 --- Сиб. Матем. Жур., 2000, Т.~41, \No~1, 150-163.

\item\label{2001-2}
{\em   Линке Ю.Ю. Саханенко А.И.}  Явное асимптотически нормальное
оценивание  параметров уравнения Михаэлиса--Ментен. --- Сиб. Матем. Жур., 2001, Т.~42, \No~3, 610-633.

\item\label{2009}
{\em Линке Ю.Ю., Саханенко А.И.} Асимптотически оптимальное
оценивание  в задаче линейной регрессии при невыполнении некоторых
классических предположений. --- Сиб. Матем. Жур.,
2009,   Т.~50, \No~2, 380-396.

\item\label{2008}
 {\em Линке Ю.Ю. Саханенко А.И.}  Асимптотически нормальное оценивание
параметра в задаче дробно -- линейной регрессии со случайными
ошибками в коэффициентах.  --- Сиб. Матем. Жур., 2008,
Т.~49, \No~3, 592-619.

\item\label{1980-Seb}  {\em  Себер Дж. }  Линейный регрессионный анализ. М.: Мир, 1980.

\item\label{2014-S}  {\em  Савинкина Е.М., Саханенко А.И. } Явные оценки неизвестного параметра в одной задаче степенной регрессии. ---
 Сиб. электрон. матем. изв., 2014, Т.11, 725-733.

\item\label{1971-F}   {\em Федоров В.В.}  Теория оптимального эксперимента. М.: Наука, 1971.

\item\label{1984-H} {\em Хьюбер П.} Робастность в статистике.  М.:Мир, 1984.

%
%
%

  \item\label{2011-B}   {\em Bergesioa A.,   Yohaia V.} Projection Estimators for Generalized Linear Models. --- J. Amer. Stat. Assoc., 2011,   V. 106(494), 661-671.

\item\label{1975-B} {\em Bickel P.J.}  One-step Huber Estimates in the Linear Model. ---
J. Amer. Stat. Assoc. 1975, V.~70, 428-434.
%
%
%
%
%
%

\item\label{1982-C} {\em
 Carrol  R.J., Ruppert D.} Robust estimation in heteroscedastic
linear models. --- Ann.  Statist.,  1982, V.~10,  2,  429-441.

%
%
%

\item\label{1999-F1}{\em Fan J., and Chen, J.} One-step local quasi-likelihood estimation. --- Journal of Royal Statistical Society B, 1999, 61, 927-943.

\item\label{1999-F} {\em Fan, J. and Jiang, J.}  Variable bandwidth and one-step local M-estimator. --- Science in China, 1999, (Series A), 29, 1-15.

%
%
%

\item\label{1925-F} {\em Fisher R.A. } Theory of statistical estimation. ---   Proc. Camb.
Phil.Soc., 1925, Soc.22,  700-725.

%

\item\label{1997-H} {\em Heyde C.C.} Quasi-likelihood and its application: a general approach to optimal parameter estimation. 1997. Springer.

%
%


%

%
%

\item\label{1985-J} {\em Janssen P., Jureckova J., Veraverbeke N.}  Rate of convergence
of one- and two-step M-estimators with applications to maksimum
likelihood and Pitman estimators. ---  Ann. Stat., 1985,  V.~13, N~3,
 1222-1229.
%

%

\item\label{1987-J} {\em Jureckova J.,  Portnoy S.} Asymptotics for one-step M-estimators in regression
with application to combining efficiency and high breakdown point. --- Comm. Statist. Theory Methods, 1987, 16, 2187-2200.

\item\label{1990-J} {\em Jureckova J., Sen P.K.}  Effect of the initial estimator on the
asymptotic behavior of the one-step M-estimator. ---  Ann. Inst.
Statist. Math. 1990, 42(2),  345-357.

\item\label{2006-J} {\em Jureckova J.,  Picek J. } Robust statistical methods with R.
Chapman and Hall, 2006.

\item\label{2012-J1} {\em  Jureckova J.,  Sen P. K.,  Picek J.}
Methodology in Robust and Nonparametric Statistics. Chapman and Hall, Boca Raton, London 2012.

\item\label{2012-J2} {\em Jureckova J.} Tail-behavior of estimators and of their one-step
versions. --- Journal de la Societe Francaise de Statistique, 2012, V. 153  (1), 44-51.
%

%
%
%
%
%

\item\label{1956-L} {\em Le Cam L. } On the asymptotic theory of estimation and testing hypotheses. ---
 Proceedings of the Third Berkeley Symposium on mathematical statistics and probability, 1956.



\item\label{2011-M} {\em
 Michaelis L., Menten M. L., Johnson K. A., Goody R. S.} The original Michaelis constant: translation of the 1913 Michaelis-Menten paper. ---  Biochemistry.  2011.  V. 50(39),  8264-8269.

\item\label{1994-M1} {\em Muller, Ch.H.}  One-step-M-estimators in conditionally contaminated linear models. ---  Stat. Decis., 1994,  12, 331-342.
%
%
%
\item\label{1985-R} {\em Reeds J.}  Asymptotic numbers of roots of Cauchy location likelihood equation. --- Ann. Statist., 1985, 13, 775-784.

%

%
%
%


%

\item\label{2003-S} {\em Seber G.A.F., Wild C.J.} Nonlinear Regression. 2003. John Wiley
and Sons.

\item\label{1980-S} {\em Serfling R.J.} Approximation theorems of mathematical statistics. John Willey and Song,  1980.

\item\label{1992-S} {\em Simpson, D. G., Ruppert, D. and Carroll, R. J.} On One-Step GM Estimates
and Stability of Inferences in Linear Regression. --- J. Amer. Stat. Assoc., 1992, 87, 439-450.

%
%
%

\item\label{2007-V} {\em Verrill S.}  Rate of convergence of k-step Newton estimators to
efficient likelihood estimators. ---  Statistics and Probability
Letters, 2007, 77, 1371--1376.

%

\item\label{2002-W} {\em  Welsh A.H., Ronchetti E.} A journey in single steps: robust one-step M-estimation in linear regression. ---
Journal of Statistical Planning and Inference, 2002, V. 103  (1-2),  287-310.

\end{enumerate}

\end{document}